\documentclass[a4paper,12pt]{amsart}
\usepackage{amsmath,amsthm,amssymb,latexsym,epic}
\usepackage{graphicx,enumerate}
\usepackage{comment}

\newtheorem{theorem}{Theorem}
\newtheorem{lemma}[theorem]{Lemma}
\newtheorem{remark}[theorem]{Remark}
\newtheorem{corollary}[theorem]{Corollary}
\newtheorem{proposition}[theorem]{Proposition}
\newtheorem{example}[theorem]{Example}

\newtheorem{conjecture}[theorem]{Conjecture}
\usepackage[all]{xy}

\newcommand{\tto}{\twoheadrightarrow}

\begin{document}
\title[Blocks and modules for Whittaker pairs]{Blocks 
and modules for Whittaker pairs}
\author{Punita Batra and Volodymyr Mazorchuk}
\date{\today}

\maketitle

\begin{abstract}
Inspired by recent activities on Whittaker modules over
various (Lie) algebras we describe some general framework 
for the study of Lie algebra modules locally finite over a 
subalgebra. As a special case we obtain a very general
setup for the study of Whittaker modules, which includes,
in particular, Lie algebras with triangular decomposition and
simple Lie algebras of Cartan type. We describe some
basic properties of Whittaker modules, including a block 
decomposition of the category of Whittaker modules and 
certain properties of simple Whittaker modules under some 
rather mild assumptions. We establish a connection between 
our general setup and the general setup of Harish-Chandra 
subalgebras in the sense of Drozd, Futorny and Ovsienko. 
For Lie algebras with triangular decomposition we construct 
a family of simple Whittaker modules (roughly depending on
the choice of a pair of weights in the dual of the Cartan 
subalgebra), describe their annihilators  and formulate 
several classification conjectures. In particular, we construct 
some new simple Whittaker modules for the Virasoro algebra. 
Finally, we construct a series of simple Whittaker modules 
for the Lie algebra of derivations of the polynomial algebra, 
and consider several finite dimensional examples, where we study 
the category of  Whittaker modules over solvable Lie algebras 
and  their relation to  Koszul algebras.
\end{abstract}

\section{Introduction}\label{s1}

The original motivation for this paper stems from the recent activities
on Whittaker modules for some infinite dimensional (Lie) algebras, which
resulted in the papers \cite{On1,On2,Ta,BO,Ch,OW,TZ,LW,Wa,LWZ,WZ}.
Whittaker modules for the Lie algebra $\mathfrak{sl}_2$ appear first
in the work \cite{AP} of Arnal and Pinczon. For all simple 
finite dimensional complex Lie algebras they were constructed 
by Kostant in \cite{Ko}. As in the ``classical'' highest weight
representation theory of Lie algebras, Whittaker modules are associated 
to a fixed triangular decomposition $\mathfrak{g}=\mathfrak{n}_-\oplus
\mathfrak{h}\oplus\mathfrak{n}_+$. However, in contrast to this 
``classical'' theory, Whittaker modules are not weight modules, 
which means that the action of nonzero elements from the Cartan 
subalgebra $\mathfrak{h}$ on Whittaker modules is (usually) not 
diagonalizable. For simple finite dimensional complex Lie algebras
Whittaker modules subsequently appeared in connection to parabolic
induction and related generalizations of the BGG category $\mathcal{O}$,
see \cite{MD1,MD2,Ba,MS1,MS2,KM}.

In the recent papers mentioned above Whittaker modules are constructed
and studied for some deformations of $\mathfrak{sl}_2$, in particular,
for generalized Weyl algebras; and for several infinite dimensional 
Lie algebras, in  particular, for the Virasoro algebra and its  
various generalizations. At the same time, there are many other Lie 
algebras for which one can naturally define Whittaker modules, notably,
affine Kac-Moody algebras, Witt algebras, Lie algebras of Cartan
type. Our main motivation for the present paper was an attempt to
understand the general picture of ``Whittaker type modules'' for
Lie algebras. We managed to find only one attempt to define some 
general setup in \cite{Wa}, however, the setup of that paper is still
rather restrictive and mostly directed to some generalizations of the
Virasoro algebra (for example, it does not cover the general case
of affine Kac-Moody algebras and even simple finite dimensional Lie
algebras). 

Comparison of the methods used in the papers listed above suggests the
following common feature of the situation: we have some Lie algebra
$\mathfrak{g}$ (possibly infinite dimensional) and some subalgebra
$\mathfrak{n}$ of $\mathfrak{g}$ (also possibly infinite dimensional,
and usually a kind of ``nilpotent''). Whittaker modules associated to
the pair $\mathfrak{n}\subset \mathfrak{g}$ are $\mathfrak{g}$-modules,
generated by a one-dimensional $\mathfrak{n}$-invariant subspace.
Moreover, the pair $\mathfrak{n}\subset \mathfrak{g}$ is chosen such
that the action of $\mathfrak{n}$ on any Whittaker $\mathfrak{g}$-module
is locally finite. Generalizing this situation it is natural to
start with a  pair $\mathfrak{n}\subset \mathfrak{g}$ of Lie algebras
and look for $\mathfrak{g}$-modules, the action of $\mathfrak{n}$ on which
is locally finite. This is our general setup in the present paper. Apart from
Whittaker modules this setup also includes classical weight modules
(in the case when the subalgebra $\mathfrak{n}$ is chosen to be a Cartan
subalgebra of $\mathfrak{g}$), and classical Harish-Chandra and
generalized Harish-Chandra modules, see e.g. \cite{Vo,PZ}. However, 
we will see that some of the results, repeatedly appearing in the papers
devoted to Whittaker modules in special cases, may be proved already
in such a general setup, or under some mild, but still general restrictions.

Let us now briefly describe the content of the paper. In Section~\ref{s2}
we describe the general setup to study $\mathfrak{g}$-modules locally 
finite over a subalgebra $\mathfrak{n}$. We show that under some 
restrictions such modules form a Serre subcategory in the category of 
all modules; and that examples of such modules can be obtained using the 
usual induction functor provided that the adjoint module
$\mathfrak{g}/\mathfrak{n}$ is itself locally finite over $\mathfrak{n}$.
We also describe a block decomposition for the category of such modules.
In Section~\ref{s3} we study, as a special case, a general setup
for Whittaker modules, where we assume that the subalgebra
$\mathfrak{n}$ is quasi-nilpotent and acts locally nilpotent on the
adjoint module $\mathfrak{g}/\mathfrak{n}$. We define Whittaker modules
in this situation and describe their basic properties. 
Using the general block decomposition theorem, obtained in 
Section~\ref{s2}, we prove that blocks of the category of 
Whittaker modules trivialize in the sense that they bijectively
correspond to characters of $\mathfrak{n}$. We also show that 
in the case when $\mathfrak{n}$ is finite dimensional, the algebra
$U(\mathfrak{n})$ is a Harish-Chandra subalgebra of 
$U(\mathfrak{g})$ in the sense of Drozd, Futorny and Ovsienko
(\cite{DOF}). In Section~\ref{s4} we study simple
Whittaker modules and their connection to one-dimensional
$\mathfrak{n}$-modules. For Lie algebras with triangular
decompositions we construct a class of  simple 
Whittaker modules as submodules in completions of usual
Verma modules. We show that such simple Whittaker modules inherit 
the annihilator from the corresponding Verma module. We also
formulate some conjectures on the structure of simple
and induced  Whittaker modules. In Section~\ref{s5}
we consider a concrete example of the Lie algebra $\mathfrak{w}_n$
of derivations of a polynomial algebra. We construct two Whittaker 
pairs associated with this algebra, one with a finite dimensional 
nilpotent subalgebra and the other, a kind of the ``opposite one'', 
with an infinite dimensional quasi-nilpotent subalgebra. We briefly 
review the theory of highest weight and lowest weight modules
and use it to construct a series of simple Whittaker modules
for $\mathfrak{w}_n$. Generalizing some of the arguments we prove
a general theorem on the existence of simple Whittaker modules
with a one-dimensional space of Whittaker vectors.
Finally, we complete the paper with several
examples for solvable finite dimensional algebras in Section~\ref{s6}.
Here the description of simple modules seems to be relatively easy, so 
the most interesting question is about the category of Whittaker modules 
(say, those of finite length). We believe that blocks of this category
are in some cases described by Koszul algebras.
\vspace{5mm}

\noindent
{\bf Acknowledgments.} The work was done during the visit of the
first author to Uppsala University. The financial support and
hospitality of Uppsala University are gratefully acknowledged. For the
second author the work was partially supported by the Royal Swedish
Academy of Sciences and the Swedish Research Council.

\section{General setup}\label{s2}

\subsection{Categories of locally finite dimensional modules}\label{s2.1}

In this paper we will work over the field $\mathbb{C}$ of complex 
numbers, so all unspecified vector spaces, homomorphisms and tensor
products are taken over $\mathbb{C}$. For a Lie algebra  $\mathfrak{a}$
we denote by $U(\mathfrak{a})$ the universal enveloping algebra of
$\mathfrak{a}$. For simplicity we assume that all vector spaces 
(in particular all algebras) have at most countable dimension.
For an abelian category $\mathfrak{A}$ of modules over some algebra
we denote by $\overline{\mathfrak{A}}$ the full subcategory of 
$\mathfrak{A}$, consisting of all modules of finite length.

Let $\mathfrak{g}$ be a nonzero complex Lie algebra (possibly 
infinite dimensional) and $\mathfrak{n}$ be a subalgebra of 
$\mathfrak{g}$. If $M$ is a $\mathfrak{g}$-module, then we say that 
the action of $\mathfrak{n}$ on $M$ is {\em locally finite} provided 
that $U(\mathfrak{n})v$ is finite dimensional for all $v\in M$. 
Let $\mathfrak{W}^{\mathfrak{g}}_{\mathfrak{n}}$ denote
the full subcategory of the category $\mathfrak{g}\text{-}\mathrm{Mod}$ 
of all $\mathfrak{g}$-modules, which consists of all 
$\mathfrak{g}$-modules, the action of $\mathfrak{n}$ on which is 
locally finite. Obviously, $\mathfrak{W}^{\mathfrak{g}}_{\mathfrak{n}}$ 
is an abelian subcategory  of $\mathfrak{g}\text{-}\mathrm{Mod}$ with 
usual kernels and cokernels.

\begin{proposition}\label{prop1}
If $\mathfrak{n}$ is finite dimensional, then
$\mathfrak{W}^{\mathfrak{g}}_{\mathfrak{n}}$ is a Serre subcategory of
$\mathfrak{g}\text{-}\mathrm{Mod}$.
\end{proposition}

\begin{proof}
The nontrivial part of this claim is to show that for any
short exact sequence
\begin{displaymath}
0\to X\to Y\to Z\to 0 
\end{displaymath}
in $\mathfrak{g}\text{-}\mathrm{Mod}$ we have that 
$X,Z\in \mathfrak{W}^{\mathfrak{g}}_{\mathfrak{n}}$ implies 
$Y\in \mathfrak{W}^{\mathfrak{g}}_{\mathfrak{n}}$. Let $v\in Y$. Then
$U(\mathfrak{n})v/(U(\mathfrak{n})v\cap X)$ is finite dimensional
as $Z\in \mathfrak{W}^{\mathfrak{g}}_{\mathfrak{n}}$. 
At the same time, the set
\begin{displaymath}
I=\{u\in U(\mathfrak{n}):uv\in X\} 
\end{displaymath}
is a left ideal of $\mathfrak{n}$. As $\mathfrak{n}$ is finite dimensional,
the algebra $U(\mathfrak{n})$ is left noetherian (\cite[2.3.8]{Di}) and 
hence $I$ is finitely generated, say by some elements $u_1,\dots,u_k$. This
means that $U(\mathfrak{n})v\cap X$ coincides with the sum of 
the subspaces $U(\mathfrak{n})u_iv$, $i=1,\dots,k$, each of which is
finite dimensional as $X\in \mathfrak{W}^{\mathfrak{g}}_{\mathfrak{n}}$. 
Hence $U(\mathfrak{n})v\cap X$ and thus also $U(\mathfrak{n})v$ are 
finite dimensional as well and the claim follows.
\end{proof}

From the definition of $\mathfrak{W}^{\mathfrak{g}}_{\mathfrak{n}}$
we have that a $\mathfrak{g}$-module  $M$ belongs to 
$\mathfrak{W}^{\mathfrak{g}}_{\mathfrak{n}}$ if and only
if there is a filtration, 
\begin{equation}\label{eq1}
0=M_0\subset M_1\subset M_2\subset \dots,\quad\quad
M=\bigcup_{i=0}^{\infty}M_i,
\end{equation}
of $M$ by finite dimensional $\mathfrak{n}$-invariant subspaces.
If $X$ is a finite dimensional $\mathfrak{n}$-module and 
$L$ is a simple finite dimensional $\mathfrak{n}$-module, we denote
by $[X:L]$ the multiplicity of $L$ in $X$ (which is obviously 
well-defined). For $M$ and $L$ as above we set
\begin{displaymath}
[M:L]:=\sup_{i\geq 0}\{[M_i:L]\}\in\{\infty,0,1,2\dots\},
\end{displaymath}
and it is easy to see that $[M:L]$ does not depend on the
choice of the filtration \eqref{eq1}. Using this notation we define the
following full subcategories of $\mathfrak{W}^{\mathfrak{g}}_{\mathfrak{n}}$:
\begin{itemize}
\item the subcategory $\mathfrak{Q}^{\mathfrak{g}}_{\mathfrak{n}}$, 
consisting of all $M$ which are semisimple over $\mathfrak{n}$;
\item the subcategory $\mathfrak{G}^{\mathfrak{g}}_{\mathfrak{n}}$, 
consisting of all $M$ such that $[M:L]<\infty$ for all $L$;
\item the subcategory $\mathfrak{H}^{\mathfrak{g}}_{\mathfrak{n}}$
of {\em Harish-Chandra} modules, consisting of all $M$ which are semisimple
over $\mathfrak{n}$ and such that $[M:L]<\infty$ for all $L$.
\end{itemize}

If $\mathfrak{C}$ is any of the categories 
$\mathfrak{W}^{\mathfrak{g}}_{\mathfrak{n}}$,
$\mathfrak{Q}^{\mathfrak{g}}_{\mathfrak{n}}$,
$\mathfrak{G}^{\mathfrak{g}}_{\mathfrak{n}}$,
$\mathfrak{H}^{\mathfrak{g}}_{\mathfrak{n}}$, then it is natural 
to ask what are simple objects of $\mathfrak{C}$. One might expect that
under some natural assumptions either $\mathfrak{C}$ or 
$\overline{\mathfrak{C}}$ has a block decomposition
with possibly finitely many simple objects in each block. In the
latter case it is known that the category $\overline{\mathfrak{C}}$
can be described as the category of finite dimensional modules over some 
(completed) algebra (\cite{Ga}), and so it is natural to ask what 
this algebra is.

\begin{example}\label{ex2}
{\rm 
Let $\mathfrak{g}$ be a semisimple finite dimensional complex
Lie algebra and $\mathfrak{n}=\mathfrak{h}$ be a Cartan subalgebra of 
$\mathfrak{g}$. Then the category 
$\mathfrak{Q}^{\mathfrak{g}}_{\mathfrak{h}}$ is the category of 
all $\mathfrak{g}$-modules, which are weight with respect to 
$\mathfrak{h}$. The category $\mathfrak{H}^{\mathfrak{g}}_{\mathfrak{h}}$ 
is the category of  all weight $\mathfrak{g}$-modules with finite
dimensional weight spaces. All simple objects in 
$\mathfrak{H}^{\mathfrak{g}}_{\mathfrak{h}}$ are classified
(\cite{Ma}) and it is known that they coincide with simple objects in 
the category $\mathfrak{G}^{\mathfrak{g}}_{\mathfrak{h}}$. 
}
\end{example}

\begin{example}\label{ex3}
{\rm 
If $\mathfrak{g}$ is the Lie algebra of an algebraic group 
and $\mathfrak{n}$ is a symmetrizing Lie subalgebra of $\mathfrak{g}$, 
then $\mathfrak{H}^{\mathfrak{g}}_{\mathfrak{h}}$ consists of usual
{\em Harish-Chandra modules} (e.g. in the sense of \cite[Chapter~9]{Di}).
}
\end{example}

Some important properties of the category 
$\mathfrak{W}^{\mathfrak{g}}_{\mathfrak{h}}$ are given by the following
statements:

\begin{proposition}\label{prop4}
Let $\mathfrak{g}$ be a Lie algebra and $\mathfrak{n}$ be a 
subalgebra of $\mathfrak{g}$. Then the category 
$\mathfrak{W}^{\mathfrak{g}}_{\mathfrak{h}}$ is stable under the usual
tensor product of $\mathfrak{g}$-modules, in particular, the category
$\mathfrak{W}^{\mathfrak{g}}_{\mathfrak{h}}$ is a monoidal category.
\end{proposition}

\begin{proof}
Let $M,N\in \mathfrak{W}^{\mathfrak{g}}_{\mathfrak{h}}$, 
$v\in M$ and $w\in N$. Then 
\begin{displaymath}
U(\mathfrak{g})(v\otimes w)\subset
(U(\mathfrak{g})v)\otimes (U(\mathfrak{g})w), 
\end{displaymath}
and the latter space is finite dimensional as both
$U(\mathfrak{g})v$ and $U(\mathfrak{g})w$ are
(because $M,N\in \mathfrak{W}^{\mathfrak{g}}_{\mathfrak{h}}$). 
The claim follows.
\end{proof}

\begin{proposition}\label{prop5}
Let $\mathfrak{g}$ be a Lie algebra and $\mathfrak{n}$ be a 
subalgebra of $\mathfrak{g}$. Then every nonzero module
$M\in \mathfrak{W}^{\mathfrak{g}}_{\mathfrak{h}}$ has a 
finite dimensional simple $\mathfrak{n}$-submodule.
\end{proposition}

\begin{proof}
Take any nonzero $v\in M$. Then $U(\mathfrak{n})v$ is a nonzero 
finite dimensional $\mathfrak{n}$-submodule of $M$ and hence has a 
well-defined socle (as an $\mathfrak{n}$-module). The claim follows.
\end{proof}

\subsection{The induction functor}\label{s2.2}
Let $\mathfrak{g}$ be a Lie algebra and $\mathfrak{n}$ be a 
subalgebra of $\mathfrak{g}$. Then we have the usual restriction functor
\begin{displaymath}
\mathrm{Res}^{\mathfrak{g}}_{\mathfrak{n}}:
\mathfrak{W}^{\mathfrak{g}}_{\mathfrak{n}}\to 
\mathfrak{W}^{\mathfrak{n}}_{\mathfrak{n}},
\end{displaymath}
which is obviously exact and hence potentially might have a left adjoint
and a right adjoint. For the category of all modules, the corresponding left 
adjoint is the usual induction functor
$\mathrm{Ind}^{\mathfrak{g}}_{\mathfrak{n}}:=
U(\mathfrak{g})\otimes_{U(\mathfrak{n})}{}_-$ and the corresponding right 
adjoint is the usual coinduction functor
$\mathrm{Coind}^{\mathfrak{g}}_{\mathfrak{n}}:=
\mathrm{Hom}_{U(\mathfrak{n})}(U(\mathfrak{g}),{}_-)$. The problem is that
these functors do not have to map $\mathfrak{W}^{\mathfrak{n}}_{\mathfrak{n}}$
to $\mathfrak{W}^{\mathfrak{g}}_{\mathfrak{n}}$ in the general case. 
However, we can state at least the following:

\begin{theorem}\label{thm6}
If the adjoint $\mathfrak{n}$-module
$\mathfrak{g}/\mathfrak{n}$ belongs to 
$\mathfrak{W}^{\mathfrak{n}}_{\mathfrak{n}}$, then 
$\mathrm{Ind}^{\mathfrak{g}}_{\mathfrak{n}}$ maps
$\mathfrak{W}^{\mathfrak{n}}_{\mathfrak{n}}$ to
$\mathfrak{W}^{\mathfrak{g}}_{\mathfrak{n}}$. In particular
$(\mathrm{Ind}^{\mathfrak{g}}_{\mathfrak{n}},
\mathrm{Res}^{\mathfrak{g}}_{\mathfrak{n}})$ is an adjoint
pair of functors between $\mathfrak{W}^{\mathfrak{n}}_{\mathfrak{n}}$
and $\mathfrak{W}^{\mathfrak{g}}_{\mathfrak{n}}$.
\end{theorem}

\begin{proof}
We choose a filtration of $\mathfrak{g}/\mathfrak{n}$ of the form
\eqref{eq1}:
\begin{displaymath}
0=X_0\subset X_1\subset X_2\subset \dots, \quad\quad\quad
\bigcup_{i=0}^{\infty} X_i= \mathfrak{g}/\mathfrak{n}.
\end{displaymath}
Note that each $X_i$ is a finite dimensional $\mathfrak{n}$-module.
Now we choose a special PBW basis of $\mathfrak{g}$ (indexed
by a well-ordered at most countable set). If
$\mathfrak{n}$ is finite dimensional, we start by choosing a
basis of $\mathfrak{n}$, then extend it by some elements which
induce a basis of $X_1$, then extend the result  by some elements 
which induce a basis of $X_2/X_1$ and so on. If
$\mathfrak{n}$ is infinite dimensional, we fix some basis of 
$\mathfrak{n}$ and then alternate the elements from this basis
first with some elements which induce a basis of $X_1$, then 
with some elements which induce a basis of $X_2/X_1$ and so on.  

Let $M\in \mathfrak{W}^{\mathfrak{n}}_{\mathfrak{n}}$. Then any
element from $U(\mathfrak{g})\otimes_{U(\mathfrak{n})}M$ can be
written as a finite sum
\begin{displaymath}
x=\sum_i u_i\otimes v_i 
\end{displaymath}
for some $u_i\in U(\mathfrak{g})$ and $v_i\in M$. As this sum is finite, 
to prove that $\dim U(\mathfrak{n})x<\infty$ it is enough to prove that
$\dim U(\mathfrak{n})u_i\otimes v_i<\infty$ for any $i$. We can write 
$u_i$ as a finite linear combination of standard monomials in the PBW basis 
chosen in the previous paragraph. As this linear combination is finite,
from our choice of the basis it follows that there exists $j$ such that 
all basis elements occurring in this linear combination either belong
to $\mathfrak{n}$ or descend to elements from $X_j$. Consider now the element
$u(u_i\otimes v_i)$ for some $u\in U(\mathfrak{n})$. Commuting all 
elements from $\mathfrak{n}$ to the right and moving them through the
tensor product we obtain that $u(u_i\otimes v_i)$ lies in the vector space,
which can be identified with $X_j\otimes U(\mathfrak{n})v_i$. 
The latter is a finite dimensional vector space as 
$X_j$ is finite dimensional and 
$M\in \mathfrak{W}^{\mathfrak{n}}_{\mathfrak{n}}$.
This implies that the action of $U(\mathfrak{n})$ on 
$U(\mathfrak{g})\otimes_{U(\mathfrak{n})}M$ is locally finite.
Therefore $\mathrm{Ind}^{\mathfrak{g}}_{\mathfrak{n}}$ maps
$\mathfrak{W}^{\mathfrak{n}}_{\mathfrak{n}}$
to $\mathfrak{W}^{\mathfrak{g}}_{\mathfrak{n}}$ and
the claim of the theorem follows follows from the fact that 
$(\mathrm{Ind}^{\mathfrak{g}}_{\mathfrak{n}},
\mathrm{Res}^{\mathfrak{g}}_{\mathfrak{n}})$ is an adjoint
pair of functors between $\mathfrak{n}\text{-}\mathrm{Mod}$
and $\mathfrak{g}\text{-}\mathrm{Mod}$.
\end{proof}

\begin{corollary}\label{cor7}
Assume that $\mathfrak{g}/\mathfrak{n}\in
\mathfrak{W}^{\mathfrak{n}}_{\mathfrak{n}}$. Then for any 
finite dimensional $\mathfrak{n}$-module $V$ every submodule
of the module $\mathrm{Ind}^{\mathfrak{g}}_{\mathfrak{n}} V$
has a simple finite dimensional $\mathfrak{n}$-submodule.
\end{corollary}

\begin{proof}
This follows immediately from Theorem~\ref{thm6} and 
Proposition~\ref{prop5}.
\end{proof}

Note that, if $\mathfrak{g}$ is finite dimensional, then the condition 
of  Theorem~\ref{thm6} is obviously satisfied for any $\mathfrak{n}$,
so one gets that the property to be locally finite dimensional
with respect to the action of $U(\mathfrak{n})$ is preserved
under induction.

If $\mathfrak{n}=\mathfrak{h}$ is a Cartan subalgebra of
$\mathfrak{g}$ in the general sense (for example, in the case when
$\mathfrak{g}$ is semisimple finite dimensional, or affine Kac-Moody, or
the Virasoro algebra, or a Witt algebra, or an algebra of Cartan type),
Theorem~\ref{thm6} implies that the $\mathfrak{g}$-module, induced from a 
weight $\mathfrak{h}$-module, is a generalized weight module. In these
cases even a stronger statement is true, namely that the module induced
from a weight module is a weight module.

\subsection{Block decomposition of 
$\mathfrak{W}^{\mathfrak{g}}_{\mathfrak{n}}$}\label{s2.3}

In this subsection we assume that $\mathfrak{g}$ is a Lie algebra and
$\mathfrak{n}$ is a subalgebra of $\mathfrak{g}$ such that the adjoint $\mathfrak{n}$-module $\mathfrak{g}/\mathfrak{n}$ belongs to 
$\mathfrak{W}^{\mathfrak{n}}_{\mathfrak{n}}$. Let 
$\mathrm{Irr}_{\mathfrak{n}}^{f}$ denote the set of isomorphism
classes of simple finite dimensional $\mathfrak{n}$-modules. If
$X\subset \mathrm{Irr}_{\mathfrak{n}}^{f}$ and $L$ is a 
simple finite dimensional $\mathfrak{n}$-module, we will loosely
write $L\in X$ if $X$ contains the isomorphism class of $L$. Define
an equivalence relation, $\sim$, on $\mathrm{Irr}_{\mathfrak{n}}^{f}$
as the smallest equivalence relation satisfying the following two
conditions:
\begin{enumerate}[(I)]
\item\label{cond1} For $L,S\in \mathrm{Irr}_{\mathfrak{n}}^{f}$ we have
$L\sim S$ if there exists an indecomposable finite dimensional 
$\mathfrak{n}$-module $M$ such that both $[M:L]\neq 0$ and $[M:S]\neq 0$.
\item\label{cond2} For $L,S\in \mathrm{Irr}_{\mathfrak{n}}^{f}$ we have
$L\sim S$ if $[\mathfrak{g}/\mathfrak{n}\otimes L:S]\neq 0$.
\end{enumerate}

\begin{example}\label{ex10}
{\rm  
Let $\mathfrak{g}$ be a simple finite dimensional complex Lie algebra
and $\mathfrak{n}=\mathfrak{h}$ be a Cartan subalgebra in $\mathfrak{g}$.
Then $\mathfrak{h}$ is a commutative Lie algebra and hence 
$\mathrm{Irr}_{\mathfrak{n}}^{f}$ can be identified with the dual 
space $\mathfrak{h}^*$ in the natural way. Furthermore, if $M$ is 
an indecomposable finite dimensional $\mathfrak{h}$-module, then,
because of the commutativity of $\mathfrak{h}$, we have a decomposition
\begin{displaymath}
M=\oplus_{\lambda\in \mathfrak{h}^*} M_{\lambda}, \quad
M_{\lambda}=\{v\in M: (h-\lambda(h))^kv=0\text{ for all }
h\in\mathfrak{h}, k\gg 0\}
\end{displaymath}
into a direct sum of $\mathfrak{h}$-modules. This means that the 
condition \eqref{cond1} only says $\lambda\sim \lambda$. The condition
\eqref{cond2} gives $\lambda\sim \lambda+\alpha$ for any root $\alpha$
of $\mathfrak{g}$ with respect to $\mathfrak{h}$. It follows that the
equivalence class of $\mathfrak{h}^*$ with respect to $\sim$ has the form
$\lambda+\mathbb{Z}\Delta$, where $\Delta$ is the root system of
$\mathfrak{g}$ with respect to $\mathfrak{h}$.
}
\end{example}

For $I\in \mathrm{Irr}_{\mathfrak{n}}^{f}/\sim$ and 
$M\in  \mathfrak{W}^{\mathfrak{g}}_{\mathfrak{n}}$ denote by
$M_I$ the trace (i.e. the sum of all images) in $M$ of all modules 
of the form $\mathrm{Ind}^{\mathfrak{g}}_{\mathfrak{n}}N$, where $N$ is a 
finite dimensional $\mathfrak{n}$-module such that for any simple
finite dimensional $\mathfrak{n}$-module $L$ we have that
$[N:L]\neq 0$ implies $L\in I$.

\begin{theorem}\label{thm11}
Let $M\in  \mathfrak{W}^{\mathfrak{g}}_{\mathfrak{n}}$.
\begin{enumerate}[(i)]
\item\label{thm11.1}
For any $I\in \mathrm{Irr}_{\mathfrak{n}}^{f}/\sim$ the vector space
$M_I$ is a $\mathfrak{g}$-submodule of $M$. Moreover, for any simple
finite dimensional $\mathfrak{n}$-module $L$ we have
$[M_I:L]\neq 0$ implies $L\in I$.
\item\label{thm11.2}
We have $M=\oplus_{I\in \mathrm{Irr}_{\mathfrak{n}}^{f}/\sim} M_I$.
\item\label{thm11.3}
If $I,J\in \mathrm{Irr}_{\mathfrak{n}}^{f}/\sim$ and $I\neq J$
then $\mathrm{Hom}_{\mathfrak{g}}(M_I,M_J)=0$.
\end{enumerate} 
\end{theorem}

\begin{proof}
The vector space $M_I$ is a $\mathfrak{g}$-submodule as the trace of 
some $\mathfrak{g}$-modules in any $\mathfrak{g}$-module is a 
$\mathfrak{g}$-submodule. Let $N$ be a finite dimensional 
$\mathfrak{n}$-module such that for any simple finite 
dimensional $\mathfrak{n}$-module $L$ we have that $[N:L]\neq 0$ 
implies $L\in I$. We claim that for any simple finite 
dimensional $\mathfrak{n}$-module $L$ we even 
have that $[\mathrm{Ind}^{\mathfrak{g}}_{\mathfrak{n}}N:L]\neq 0$ 
implies $L\in I$.

First we recall that 
$\mathrm{Ind}^{\mathfrak{g}}_{\mathfrak{n}}N=U(\mathfrak{g})
\otimes_{U(\mathfrak{n})}N$, so we can use the PBW Theorem.
Fix some basis in $\mathfrak{g}$  as described in the proof of 
Theorem~\ref{thm6} and for $n\in\mathbb{N}$ denote by 
$U(\mathfrak{g})_n$ the linear subspace of $U(\mathfrak{g})$,
generated by all standard monomials of degree at most $n$. 
Then $U(\mathfrak{g})_n$ is an $\mathfrak{n}$-submodule of
$U(\mathfrak{g})$ with respect to the adjoint action.
Moreover, every finite dimensional $\mathfrak{n}$-submodule
of the module $U(\mathfrak{g})\otimes_{U(\mathfrak{n})}N$ is a
submodule of the $\mathfrak{n}$-module
$U(\mathfrak{g})_n\otimes_{U(\mathfrak{n})}N$ for some $n$.
By \cite[2.4.5]{Di} the $\mathfrak{g}$-module 
$U(\mathfrak{g})_n/U(\mathfrak{g})_{n-1}$
is isomorphic to the $n$-th symmetric power of $\mathfrak{g}$,
which is a submodule of
\begin{displaymath}
\mathfrak{g}^{\otimes n}:=
\underbrace{\mathfrak{g}\otimes\mathfrak{g}\otimes\dots\otimes
\mathfrak{g}}_{n\text{ factors }}.
\end{displaymath}
This implies that any simple subquotient of the $\mathfrak{n}$-module 
\begin{displaymath}
U(\mathfrak{g})_n\otimes_{U(\mathfrak{n})}N/
U(\mathfrak{g})_{n-1}\otimes_{U(\mathfrak{n})}N
\end{displaymath}
is a subquotient of the module 
$(\mathfrak{g}/\mathfrak{n})^{\otimes n}\otimes N$.
From the definition of the relation $\sim$ it follows by induction on $n$ 
that  any simple finite dimensional $\mathfrak{n}$-subquotient $L$ 
of $(\mathfrak{g}/\mathfrak{n})^{\otimes n}\otimes N$ belongs to $I$. 
This proves the claim \eqref{thm11.1}.
The claim \eqref{thm11.3} follows directly from the definitions
and the claim \eqref{thm11.1}.

To prove the claim \eqref{thm11.2} let 
$X=\oplus_{I\in \mathrm{Irr}_{\mathfrak{n}}^{f}/\sim} M_I$.
Then $X$ is a submodule of $M$ by the claim \eqref{thm11.1}. Consider
a filtration of $M$ of the form \eqref{eq1}.
As every $M_i$ is finite dimensional, we obviously have $X\cap M_i=M_i$,
which implies $X=M$. This completes the proof of the claim \eqref{thm11.2}
and of the whole theorem.
\end{proof}

For $I\in \mathrm{Irr}_{\mathfrak{n}}^{f}/\sim$ denote by 
$\mathfrak{W}^{\mathfrak{g}}_{\mathfrak{n}}(I)$ the full subcategory
of $\mathfrak{W}^{\mathfrak{g}}_{\mathfrak{n}}$, which consists of
modules $M$ satisfying $M=M_I$. The categories 
$\mathfrak{W}^{\mathfrak{g}}_{\mathfrak{n}}(I)$ will be called
{\em blocks} of $\mathfrak{W}^{\mathfrak{g}}_{\mathfrak{n}}$ and
Theorem~\ref{thm11} gives a coproduct block decomposition of 
$\mathfrak{W}^{\mathfrak{g}}_{\mathfrak{n}}$. From the category 
theoretical point of view it is enough to study each block 
$\mathfrak{W}^{\mathfrak{g}}_{\mathfrak{n}}(I)$ separately.

\begin{example}\label{ex12}
{\rm  
Let $\mathfrak{g}$ and $\mathfrak{n}$ be as in Example~\ref{ex10}.
Then Theorem~\ref{thm11} gives the usual block decomposition for 
the category of weight modules with respect to supports
of modules, indexed by cosets from  $\mathfrak{h}^*/\mathbb{Z}\Delta$.
}
\end{example}

\subsection{$\mathfrak{n}$-socles of modules in
$\mathfrak{W}^{\mathfrak{g}}_{\mathfrak{n}}$}\label{s2.4}

Let $\mathfrak{g}$ be a Lie algebra and $\mathfrak{n}$ be
a Lie subalgebra of $\mathfrak{g}$. Then every $M\in 
\mathfrak{W}^{\mathfrak{g}}_{\mathfrak{n}}$ has a well-defined
$\mathfrak{n}$-socle $\mathrm{soc}_{\mathfrak{n}}(M)$. As usual
we have the following standard result:

\begin{proposition}\label{prop101}
Let $M\in \mathfrak{W}^{\mathfrak{g}}_{\mathfrak{n}}$ and $L$
be a simple finite dimensional $\mathfrak{n}$-module. Then 
\begin{displaymath}
[\mathrm{soc}_{\mathfrak{n}}(M):L]=\dim
\mathrm{Hom}_{\mathfrak{n}}(L,M)
=\dim
\mathrm{Hom}_{\mathfrak{g}}(\mathrm{Ind}_{\mathfrak{n}}^{\mathfrak{g}}L,M).
\end{displaymath}
\end{proposition}

\begin{proof}
This follows from the usual adjunction between induction and
restriction.
\end{proof}

Note that in the case $\mathfrak{g}/\mathfrak{n}\in
\mathfrak{W}^{\mathfrak{n}}_{\mathfrak{n}}$, the rightmost
homomorphism space in the formulation of Proposition~\ref{prop101} 
is inside the category $\mathfrak{W}^{\mathfrak{g}}_{\mathfrak{n}}$ by
Theorem~\ref{thm6}.

\section{Whittaker modules revisited}\label{s3}

\subsection{Finite dimensional representations of quasi-nilpotent
Lie algebras}\label{s3.1}

For a Lie algebra $\mathfrak{n}$ define inductively ideals
$\mathfrak{n}_0:=\mathfrak{n}$ and
$\mathfrak{n}_i:=[\mathfrak{n}_{i-1},\mathfrak{n}]$, $i>0$. Then 
we have a sequence of ideals
\begin{displaymath}
\mathfrak{n}= \mathfrak{n}_0\supset \mathfrak{n}_1\supset
\mathfrak{n}_2\supset \dots.
\end{displaymath}
We will say that $\mathfrak{n}$ is {\em quasi-nilpotent}
provided that $\displaystyle \bigcap_{i=0}^{\infty} \mathfrak{n}_i=0$.
Obviously, any nilpotent Lie algebra  is quasi-nilpotent. 
Until the end of this subsection we assume that $\mathfrak{n}$ is 
quasi-nilpotent. The next example of a quasi-nilpotent Lie algebra,
which is no nilpotent, comes from the Virasoro algebra.

\begin{example}\label{ex14}
{\rm  
Let $\mathfrak{n}$ have the basis $\{e_i:i=1,2,3,\dots\}$ with the
Lie bracket $[e_i,e_j]=(j-i)e_{i+j}$. Then it is easy to see that 
$\mathfrak{n}_k$, $k>1$, is the subspace of $\mathfrak{n}$, spanned by
$\{e_i:i=k+2,k+3,k+4,\dots\}$ and hence $\mathfrak{n}$ is quasi-nilpotent.
}
\end{example}

Example~\ref{ex14} generalizes in a straightforward way to the 
positive part $\mathfrak{n_+}$ of any Lie  algebra $\mathfrak{g}$ with a
triangular decomposition $\mathfrak{g}=\mathfrak{n}_-\oplus
\mathfrak{h}\oplus\mathfrak{n}_+$ in the sense of \cite{MP}.

The main for us property of quasi-nilpotent Lie algebras is the
following easy fact:

\begin{proposition}\label{prop15}
Let $\mathfrak{n}$ be a quasi-nilpotent Lie algebra and
$M$ be a finite dimensional $\mathfrak{n}$-module. Then
there is $i\in\mathbb{N}$ such that $\mathfrak{n}_iM=0$.
\end{proposition}

\begin{proof}
Let $\mathfrak{i}$ be the kernel of the Lie algebra homomorphism 
from $\mathfrak{n}$ to $\mathrm{End}_{\mathbb{C}}(M)$, defining 
the $\mathfrak{n}$-module  structure on $M$. Then $\mathfrak{i}$ is an 
ideal of $\mathfrak{n}$ of finite codimension. Set 
$a_i=\dim (\mathfrak{n}_i+\mathfrak{i}/\mathfrak{i})$.
As $\mathfrak{n}_{i+1}\subset \mathfrak{n}_i$, we have that the
sequence $\mathbf{a}:=\{a_i\}$ is weakly decreasing. As 
$\mathfrak{i}$ has finite codimension in $\mathfrak{n}$, all elements in
$\mathbf{a}$ are finite positive integers. As 
$\displaystyle \bigcap_{i=0}^{\infty} \mathfrak{n}_i=0$,
the sequence $\mathbf{a}$ converges to $1$. This means that there exists
$i\in\mathbb{N}$ such that $a_i=1$, which means 
$\mathfrak{n}_i\subset \mathfrak{i}$ and thus $\mathfrak{n}_iM=0$.
This completes the proof.
\end{proof}

\begin{corollary}\label{cor16}
Let $\mathfrak{n}$ be a quasi-nilpotent Lie algebra and
$L$ be a simple finite dimensional $\mathfrak{n}$-module. Then we have:
\begin{enumerate}[(i)]
\item\label{cor16.1} $\dim L=1$.
\item\label{cor16.2} $[\mathfrak{n},\mathfrak{n}]L=0$.
\end{enumerate}
\end{corollary}

\begin{proof}
By Proposition~\ref{prop15}, we have $\mathfrak{n}_iL=0$ for some
$i\in\mathbb{N}$. Hence $L$ is a simple finite dimensional module 
over the nilpotent Lie algebra  $\mathfrak{n}/\mathfrak{n}_i$.
Now the claim \eqref{cor16.1} follows from the Lie Theorem 
(\cite[1.3.13]{Di}). As $\dim L=1$, we also have that the Lie algebra 
$\mathrm{End}_{\mathbb{C}}(L)$ is commutative, which implies 
the claim \eqref{cor16.2}. This completes the proof.
\end{proof}

Set $\mathfrak{h}=\mathfrak{n}/[\mathfrak{n},\mathfrak{n}]$. Then
$\mathfrak{h}$ is a commutative Lie algebra. From Corollary~\ref{cor16} 
it follows that $\mathrm{Irr}_{\mathfrak{n}}^f$ can 
be naturally identified with $\mathfrak{h}^*$ 
(compare with Example~\ref{ex10}), which, in turn, can be identified
with Lie algebra homomorphisms from $\mathfrak{n}$ to $\mathbb{C}$.
In what follows we consider elements from $\mathfrak{h}^*$ as 
Lie algebra homomorphisms from $\mathfrak{n}$ to $\mathbb{C}$ under
this identification. Now we would like to establish block decomposition
for finite dimensional $\mathfrak{n}$-modules and show that the
blocks are indexed by elements from $\mathfrak{h}^*$ in the natural way.

Let $\mathfrak{n}$ be a quasi-nilpotent Lie algebra and
$M$ be a finite dimensional $\mathfrak{n}$-module. For
$\lambda\in \mathfrak{h}^*$ set
\begin{displaymath}
M_{\lambda}:=\{v\in M:(x-\lambda(x))^kv=0\text{ for all }
x\in\mathfrak{n}\text{ and }k\gg 0\}.
\end{displaymath}

\begin{proposition}\label{prop17}
\begin{enumerate}[(i)]
\item \label{prop18} Each $M_{\lambda}$ is a submodule of $M$.
\item \label{prop19} $M=\oplus_{\lambda\in \mathfrak{h}^*}M_{\lambda}$.
\end{enumerate}
\end{proposition}

\begin{proof}
Let $\mathfrak{m}$ denote the image of $\mathfrak{n}$ under the Lie
algebra homomorphism defining the module structure of $M$. We can 
consider the algebra $\mathfrak{m}$ instead of the algebra 
$\mathfrak{n}$. Then, by Proposition~\ref{prop15}, the algebra 
$\mathfrak{m}$ is a finite dimensional nilpotent Lie algebra
and the claim follows from \cite[1.3.19]{Di}.
\end{proof}

As an immediate corollary we have the following:

\begin{corollary}\label{cor18}
Let $\mathfrak{n}$ be a quasi-nilpotent Lie algebra,
$M$ be an indecomposable finite dimensional $\mathfrak{n}$-module,
and $L$ be a simple submodule of $M$. Then $[M:L]=\dim M$. In
particular, if $S$ is a simple $\mathfrak{n}$-module such that 
$S\not\cong L$, then $[M:S]=0$.
\end{corollary}

\begin{remark}\label{rem18}
{\rm
Corollary~\ref{cor18} means that for quasi-nilpotent Lie algebras
the condition \eqref{cond1} from the definition of the equivalence
relation $\sim$ (see Subsection~\ref{s2.3}) trivializes (i.e.
gives the equality relation).
}
\end{remark}

\subsection{General Whittaker setup}\label{s3.2}

Now we are getting closer to the general situation in which one can
consider Whittaker modules (see \cite{Ko,OW,TZ,LW,Wa,LWZ,WZ}). We define 
the general Whittaker setup as follows: consider a Lie algebra 
$\mathfrak{g}$ and a quasi-nilpotent Lie subalgebra $\mathfrak{n}$ of 
$\mathfrak{g}$ such that the action of $\mathfrak{n}$ on the adjoint
$\mathfrak{n}$-module $\mathfrak{g}/\mathfrak{n}$ is locally nilpotent
(in particular, $\mathfrak{g}/\mathfrak{n}\in 
\mathfrak{W}^{\mathfrak{n}}_{\mathfrak{n}}$).
In this case we will say that $(\mathfrak{g},\mathfrak{n})$ is
a {\em Whittaker pair}. If $(\mathfrak{g},\mathfrak{n})$ is
a Whittaker pair, then objects in 
$\mathfrak{W}^{\mathfrak{g}}_{\mathfrak{n}}$ will be called
{\em Whittaker modules}.
 
Our terminology corresponds to the one used in \cite{Ba,MS1,MS2}.
In \cite{Ko,OW} and many other papers a Whittaker module is additionally
supposed to be generated by a Whittaker vector (see Section~\ref{s4.1}).
The advantage of our definition is that the category of all Whittaker
modules is abelian (see Subsection~\ref{s2.1}).

\begin{example}\label{ex1901}
{\rm 
Let $\mathfrak{g}$ be a Lie algebra and $\mathfrak{z}$ be the center
of $\mathfrak{g}$. Then $(\mathfrak{g},\mathfrak{z})$ is a Whittaker pair.
} 
\end{example}

\begin{example}\label{ex19010}
{\rm 
Let $\mathfrak{n}$ be a quasi-nilpotent Lie algebra and 
$\mathfrak{a}$ be any Lie algebra. Then 
$(\mathfrak{a}\oplus \mathfrak{n},\mathfrak{n})$ is a Whittaker pair.
} 
\end{example}

\begin{example}\label{ex1902}
{\rm 
Let $\mathfrak{g}$ be a solvable Lie algebra and 
$\mathfrak{n}=[\mathfrak{g},\mathfrak{g}]$. Then 
$\mathfrak{n}$ is nilpotent (\cite[1.7.1]{Di}) 
and $(\mathfrak{g},\mathfrak{n})$ is a Whittaker pair.
} 
\end{example}

\begin{example}\label{ex19}
{\rm 
Let $\mathfrak{g}$ be a Lie algebra with a fixed triangular decomposition
$\mathfrak{g}=\mathfrak{n}_-\oplus\mathfrak{h}\oplus\mathfrak{n}_+$
in the sense of \cite{MP}. Let further $\mathfrak{n}=\mathfrak{n}_+$.
Then $(\mathfrak{g},\mathfrak{n})$ is a Whittaker pair. Note that in this
case the algebra $\mathfrak{g}$ as well as the subalgebra $\mathfrak{n}$
may be infinite dimensional. This example contains, as special cases, 
situations studied in the articles \cite{Ko,MD1,MD2,OW}, where some simple
Whittaker modules were described. It also includes many
new examples, e.g. where $\mathfrak{g}$ is an affine Kac-Moody algebra.
} 
\end{example}

\begin{example}\label{ex20}
{\rm 
Let $\mathfrak{g}$ be a simple finite dimensional Lie algebra
with a fixed Cartan subalgebra $\mathfrak{h}$. Then $\mathfrak{h}$ is
commutative, in particular, it is quasi-nilpotent. 
However,  $(\mathfrak{g},\mathfrak{h})$ is {\em not} a 
Whittaker pair as the adjoint action of $\mathfrak{h}$ on 
$\mathfrak{g}/\mathfrak{h}$ is not locally nilpotent.
} 
\end{example}

\begin{example}\label{ex21}
{\rm 
Let $\mathfrak{g}$ be a simple finite dimensional Lie algebra
with a fixed Cartan subalgebra $\mathfrak{h}$. Let $\alpha$ be
a root of $\mathfrak{g}$ with respect to $\mathfrak{h}$ and 
$\mathfrak{n}$ be the corresponding root subspace $\mathfrak{g}_{\alpha}$ 
of $\mathfrak{g}$. Then $\mathfrak{n}$ is an abelian subalgebra and 
$(\mathfrak{g},\mathfrak{n})$ is a Whittaker pair.
} 
\end{example}

\begin{example}\label{ex22}
{\rm 
Let $\mathfrak{g}$ be the Lie algebra with the basis 
\begin{displaymath}
\{e_i:i\in\{\dots,
-2,-1,0,1\}\} 
\end{displaymath}
and the Lie bracket given by $[e_i,e_j]=(j-i)e_{i+j}$.
Let $\mathfrak{n}=\langle e_1\rangle$. Then 
$(\mathfrak{g},\mathfrak{n})$ is a Whittaker pair.
This example will be considered in more details in Section~\ref{s5}.
} 
\end{example}

\begin{example}\label{ex23}
{\rm 
Let $\mathfrak{g}$ be the Lie algebra with the basis 
\begin{displaymath}
\{e_i:i\in\{-1,0,1,2,\dots\}\}
\end{displaymath}
and the Lie bracket given by $[e_i,e_j]=(j-i)e_{i+j}$.
Let $\mathfrak{n}=\langle e_1,e_2,\dots\rangle$ (see Example~\ref{ex14}). 
Then  $(\mathfrak{g},\mathfrak{n})$ is a Whittaker pair. Simple
Whittaker modules in this case were completely classified in 
\cite{Ru}.
} 
\end{example}

\begin{remark}\label{rem25}
{\rm  
It is easy to check that situations considered in the papers
\cite{TZ,LW,Wa,LWZ} also correspond to Whittaker pairs.
Examples~\ref{ex22} and \ref{ex23} generalize to any Witt algebra
$\mathfrak{w}_n$, that is the Lie algebra of derivations 
of $\mathbb{C}[t_1,\dots,t_n]$ and more generally, to some other 
infinite dimensional Lie algebras of Cartan type (see also \cite{Ru} 
for some partial results on Whittaker modules in these cases). The 
example of the algebra $\mathfrak{w}_n$ will 
be considered in more details in Section~\ref{s5}.
}
\end{remark}

\subsection{Blocks for Whittaker pairs}\label{s3.3}

From now on, if not explicitly stated otherwise, we assume that
$(\mathfrak{g},\mathfrak{n})$ is a Whittaker pair. Our first aim
is to show that in this case the block decomposition of 
$\mathfrak{W}^{\mathfrak{g}}_{\mathfrak{n}}$, described in 
Subsection~\ref{s2.3}, trivializes in the following sense:

\begin{theorem}\label{thm26}
Let $(\mathfrak{g},\mathfrak{n})$ be a Whittaker pair.
Then the equivalence relation $\sim$ from Subsection~\ref{s2.3}
is the equality relation.
\end{theorem}

\begin{proof}
From Remark~\ref{rem18} we already know that the condition \eqref{cond1}
trivializes as $\mathfrak{n}$ is quasi-nilpotent. Hence we only have
to check that the condition \eqref{cond2} trivializes as well.

Set $\mathfrak{h}=\mathfrak{n}/[\mathfrak{n},\mathfrak{n}]$ and
identify $\mathfrak{h}^*$ with Lie algebra homomorphisms from
$\mathfrak{n}$ to $\mathbb{C}$. For $\lambda\in \mathfrak{h}^*$
let $L_{\lambda}$ denote the simple one-dimensional
$\mathfrak{n}$-module, given by $\lambda$. Let $v_{\lambda}$ 
be some fixed nonzero element in $L_{\lambda}$.

To prove the claim we have to show that for any $u\in \mathfrak{n}$
the element $u-\lambda(u)$ acts locally nilpotent on the
$\mathfrak{n}$-module $\mathfrak{g}/\mathfrak{n}\otimes L_{\lambda}$.
For any $w\in \mathfrak{g}/\mathfrak{n}$ we have 
\begin{displaymath}
(u-\lambda(u))(w\otimes v_{\lambda})=
[u,w]\otimes v_{\lambda} + w\otimes u(v_{\lambda}) -
\lambda(u) w\otimes v_{\lambda}=
[u,w]\otimes v_{\lambda},
\end{displaymath}
as $(u-\lambda(u))v_{\lambda}=0$. This implies, by induction, that 
\begin{equation}\label{eq12321}
(u-\lambda(u))^k(w\otimes v_{\lambda})=\mathrm{ad}^k_u(w)\otimes 
v_{\lambda}.
\end{equation}
As $(\mathfrak{g},\mathfrak{n})$ is a Whittaker pair, the adjoint action
of any element from $\mathfrak{n}$ on the module $\mathfrak{g}/\mathfrak{n}$
is locally nilpotent. Hence $\mathrm{ad}^k_u(w)\otimes v_{\lambda}=0$ for all
$k\gg 0$, which implies $(u-\lambda(u))^k(w\otimes v_{\lambda})=0$ for all
$k\gg 0$ by \eqref{eq12321}. The claim of the theorem follows.
\end{proof}

Theorem~\ref{thm26} says that for any Whittaker pair 
$(\mathfrak{g},\mathfrak{n})$ blocks of the category 
$\mathfrak{W}^{\mathfrak{g}}_{\mathfrak{n}}$ of 
Whittaker modules, as defined in  Subsection~\ref{s2.3},
are indexed by $\lambda\in (\mathfrak{n}/[\mathfrak{n},\mathfrak{n}])^*$
in the natural way. We will denote these blocks by $\mathfrak{W}^{\mathfrak{g}}_{\mathfrak{n}}(\lambda)$.
After Theorem~\ref{thm26} it is natural to say that the main
problem in the theory of Whittaker modules is to describe the
categories $\mathfrak{W}^{\mathfrak{g}}_{\mathfrak{n}}(\lambda)$,
$\lambda\in (\mathfrak{n}/[\mathfrak{n},\mathfrak{n}])^*$.

From Theorem~\ref{thm26} it follows that general Whittaker setup 
which leads to  Whittaker modules is in some sense
``opposite'' to those pairs $(\mathfrak{g},\mathfrak{n})$,
for which one gets usual Harish-Chandra modules.

\subsection{Connection to Harish-Chandra subalgebras}\label{s3.5}

Recall (see \cite{DOF}) that an associative unital algebra $B$ 
is called {\em quasi-commutative} if $\mathrm{Ext}_B^1(L,S)=0$ 
for any two nonisomorphic simple finite dimensional $B$-modules 
$L$ and $S$. Let $A$ be an algebra and $B$ be a subalgebra of $A$. 
Following \cite{DOF} we say that $B$ is {\em quasi-central} 
provided that for any $a\in A$ the $B$-bimodule $BaB$ is finitely 
generated both as a left and as a right $B$-module. 
A subalgebra $B$ of $A$ is called a {\em Harish-Chandra subalgebra} 
if it is both, quasi-commutative and quasi-central.

\begin{theorem}\label{thm36}
Let $(\mathfrak{g},\mathfrak{n})$ be a Whittaker pair.  
\begin{enumerate}[(i)]
\item\label{thm36.1} The algebra $U(\mathfrak{n})$ is 
quasi-commutative.
\item\label{thm36.2} If $\dim\mathfrak{n}<\infty$, then
$U(\mathfrak{n})$ is a Harish-Chandra subalgebra of
$U(\mathfrak{g})$.
\end{enumerate}
\end{theorem}

\begin{proof}
The claim \eqref{thm36.1} follows from Proposition~\ref{prop17}.
To prove \eqref{thm36.2} we have only to show that 
$U(\mathfrak{n})$ is quasi-central. We would need the following
variation of the PBW Theorem:

\begin{lemma}\label{lem37}
Let $\mathfrak{g}$ be a Lie algebra and $\mathfrak{n}$ be a 
subalgebra of $\mathfrak{g}$. Let $\{a_i:i\in I\}$ (where $I$
is well-ordered) be some basis in $\mathfrak{n}$ and 
$\{b_j:j\in J\}$ (where $J$ is well-ordered) be a complement of 
$\{a_i\}$ to a basis of $\mathfrak{g}$. Then $U(\mathfrak{g})$ has a 
basis consisting of all elements of the form $\mathbf{b}\mathbf{a}$, 
where $\mathbf{b}$ is a standard monomial in $\{b_j\}$ and 
$\mathbf{a}$ is a standard monomial in $\{a_i\}$.
\end{lemma}

\begin{proof}
The proof is similar to the standard proof of the PBW Theorem 
and is left to the reader.
\end{proof}

\begin{lemma}\label{lem38}
Let $\mathfrak{g}$ be a Lie algebra and $\mathfrak{n}$ be a 
Lie subalgebra of $\mathfrak{g}$ such that the adjoint 
$\mathfrak{n}$-module $\mathfrak{g}/\mathfrak{n}$ belongs to 
$\mathfrak{W}_{\mathfrak{n}}^{\mathfrak{n}}$.  
Then the adjoint $\mathfrak{n}$-module 
$U(\mathfrak{g})/U(\mathfrak{n})$ belongs to 
$\mathfrak{W}_{\mathfrak{n}}^{\mathfrak{n}}$ as well.
\end{lemma}

\begin{proof}
Choose some basis in $\mathfrak{g}$ as in Theorem~\ref{thm6}
and the corresponding basis in 
$U(\mathfrak{g})$ as given by Lemma~\ref{lem37}. Then the adjoint $\mathfrak{n}$-module $U(\mathfrak{g})/U(\mathfrak{n})$
can be identified with the linear span of standard monomials
in $\{b_j\}$. Let $u\in U(\mathfrak{g})$. Then, writing $u$ in our
basis of $U(\mathfrak{g})$, we get a finite linear combination of 
standard monomials with some nonzero coefficients.
In particular, only finitely many standard monomials in  
$\{b_j\}$ show up (as factors of summands in this linear combination). 
As $\mathfrak{g}/\mathfrak{n}\in \mathfrak{W}_{\mathfrak{n}}^{\mathfrak{n}}$ 
by our assumption,
applying the adjoint action of $\mathfrak{n}$ to all these standard
monomials in $\{b_j\}$ we can produce, as summands, only finitely many
new standard monomials in $\{b_j\}$. The claim follows.
\end{proof}

For  $u\in U(\mathfrak{n})$ consider the $U(\mathfrak{n})$-bimodule 
$X=U(\mathfrak{n})uU(\mathfrak{n})$, By Lemma~\ref{lem38}, the
image of $X$ in $U(\mathfrak{g})/U(\mathfrak{n})$ is finite dimensional.
At the same time $X\cap U(\mathfrak{n})$ is an ideal of $U(\mathfrak{n})$.
As $\mathfrak{n}$ is finite dimensional, $U(\mathfrak{n})$ is noetherian
(\cite[2.3.8]{Di}). Hence $X\cap U(\mathfrak{n})$ is finitely generated
as a left $U(\mathfrak{n})$-module. This implies that $X$ is finitely 
generated as a left $U(\mathfrak{n})$-module. Applying the canonical 
antiinvolution on $\mathfrak{g}$ we obtain that $X$ is finitely 
generated as a right $U(\mathfrak{n})$-module as well.  Therefore 
$U(\mathfrak{n})$ is quasi-central  and the claim \eqref{thm36.2} 
of our theorem follows.
\end{proof}

\begin{remark}\label{rem309}
{\rm 
\begin{enumerate}[(a)]
\item\label{rem309-1}
There are natural examples of Whittaker pairs $(\mathfrak{g},\mathfrak{n})$,
where $\mathfrak{n}$ is finite dimensional while  $\mathfrak{g}$ is 
infinite dimensional, see Example~\ref{ex22} and Remark~\ref{rem25}.
\item\label{rem309-2}
For finite dimensional $\mathfrak{n}$ to prove Lemma~\ref{lem38}
one could alternatively argue using 
Propositions~\ref{prop4} and \ref{prop1} and arguments similar to
the ones used in the proof of Theorem~\ref{thm11}.
\item\label{rem309-3}
Theorem~\ref{thm6} follows from 
Lemma~\ref{lem38} and Propositions~\ref{prop4}.
\end{enumerate}
}
\end{remark}

\section{Whittaker vectors and simple Whittaker modules}\label{s4}

\subsection{Whittaker vectors, standard and 
simple Whittaker modules}\label{s4.1}

Let $(\mathfrak{g},\mathfrak{n})$ be a Whittaker pair and 
$\lambda\in (\mathfrak{n}/[\mathfrak{n},\mathfrak{n}])^*$. 
As in the previous section we denote by $L_{\lambda}$ the simple
one-dimensional $\mathfrak{n}$-module given by $\lambda$.
Let $v_{\lambda}$ be some basis element of $L_{\lambda}$.
Set $M_{\lambda}=U(\mathfrak{g})\otimes_{U(\mathfrak{n})}L_{\lambda}$
and call this module the {\em standard} Whittaker module.

Note that in \cite{Ko,OW} and some other papers the module 
$M_{\lambda}$ is called the universal Whittaker module, however
the latter wording might be slightly misleading in our setup as not every
Whittaker module is a quotient of $M_{\lambda}$ (or direct sums
of various $M_{\lambda}$'s). As an example one could take
$\mathfrak{g}$ to be a simple finite dimensional complex Lie algebra
with a fixed triangular decomposition 
$\mathfrak{g}=\mathfrak{n}_-\oplus \mathfrak{h}\oplus\mathfrak{n}_+$,
$\mathfrak{n}=\mathfrak{n}_+$, $\lambda=0$, and $M$ any module in
the BGG category $\mathcal{O}$, associated with this triangular
decomposition, which is not generated by its highest weights
(for example some projective module, which is not a Verma module).

Let $M\in \mathfrak{W}^{\mathfrak{g}}_{\mathfrak{n}}(\lambda)$.
A vector $v\in M$ is called a {\em Whittaker vector}
provided that $\langle v\rangle$ is an $\mathfrak{n}$-submodule
of $M$ (which is automatically isomorphic to $L_{\lambda}$ by
Theorem~\ref{thm26}). Obviously, all Whittaker vectors form an
$\mathfrak{n}$-submodule of $M$, which we will denote by  $W_{\lambda}(M)$.

\begin{proposition}\label{prop30}
Let $M\in \mathfrak{W}^{\mathfrak{g}}_{\mathfrak{n}}(\lambda)$.
\begin{enumerate}[(i)]
\item\label{prop30.1} $\dim W_{\lambda}(M)=\dim 
\mathrm{Hom}_{\mathfrak{g}}(M_{\lambda},M)$.
\item\label{prop30.2} If a Whittaker module $M$ contains a
unique (up to scalar) nonzero Whittaker vector, which, moreover,
generates $M$, then $M$ is a simple module.
\end{enumerate}
\end{proposition}

\begin{proof}
The claim \eqref{prop30.1} follows from Proposition~\ref{prop101}.
To prove the claim \eqref{prop30.2} let $N\subset M$ be 
a nonzero submodule. Then $N$ contains a nonzero Whittaker vector by 
Corollary~\ref{cor7}. Since such vector in $M$ is unique and
generates $M$, we have $N=M$. This implies that $M$ is simple
and proves \eqref{prop30.2}.
\end{proof}

Based on the examples from \cite{Ko,OW,Ch,TZ} it looks reasonable
to formulate the following conjectures:

\begin{conjecture}\label{conj31}
{\rm
Assume that $\mathfrak{g}$ is a Lie algebra with a fixed triangular
decomposition $\mathfrak{g}=\mathfrak{n}_-\oplus \mathfrak{h}\oplus
\mathfrak{n}_+$ in the sense of \cite{MP}. Let $L$ be a simple
Whittaker module for the Whittaker pair 
$(\mathfrak{g},\mathfrak{n}_+)$. Then $\mathrm{soc}_{\mathfrak{n}}(L)$
is a simple module.
}
\end{conjecture}

Later on we will give some evidence for Conjecture~\ref{conj31}
(in particular, in Subsection~\ref{s5.4} we show that simple
Whittaker modules with simple $\mathfrak{n}$-socle always exist). 
We note that Conjecture~\ref{conj31} does not extend to the general 
case, see  example in Subsection~\ref{s6.2}.

\begin{conjecture}\label{conj32}
{\rm
Assume that $\mathfrak{g}$ is a Lie algebra with a fixed triangular
decomposition $\mathfrak{g}=\mathfrak{n}_-\oplus \mathfrak{h}\oplus
\mathfrak{n}_+$ in the sense of \cite{MP}, $\mathfrak{n}=\mathfrak{n}_+$ 
and $\lambda\in (\mathfrak{n}/[\mathfrak{n},\mathfrak{n}])^*$. 
Then for generic $\lambda$ the center $Z(\mathfrak{g})$ of $U(\mathfrak{g})$ 
surjects onto the set of Whittaker vectors of $M_{\lambda}$ via 
$z\mapsto z\otimes v_{\lambda}$, $z\in Z(\mathfrak{g})$. 
}
\end{conjecture}

The Whittaker pair $(\mathfrak{g},\mathfrak{n})$ associated with a
fixed triangular decomposition of $\mathfrak{g}$
(see Example~\ref{ex19}, Conjecture~\ref{conj32}) seems to be
the most reasonable situation to study Whittaker modules as it looks
the most ``balanced'' one in the following sense: if the
algebra $\mathfrak{n}$ is much ``bigger'' than 
$\mathfrak{g}/\mathfrak{n}$ then standard
Whittaker modules should normally be simple (see for example 
some evidence for this in \cite{Ru}); on the other hand if 
$\mathfrak{n}$ is much ``smaller'' than $\mathfrak{g}/\mathfrak{n}$ 
then standard Whittaker modules should normally have ``too many''
simple quotients with no chance of classifying them (take for example
the situation described in Examples~\ref{ex19010} and \ref{ex21}). The main
advantage of special cases studied so far (\cite{Ko,OW,Ch,TZ}
and others) is that in those cases
the considered situation was balanced enough to give a reasonable
classification of generic simple Whittaker modules.
 
In the more general situation described in Section~\ref{s2}
already Conjecture~\ref{conj31} is not reasonable. In fact, as
mentioned in Example~\ref{ex2}, a special case of such situation is
the study of weight modules over simple complex finite dimensional
Lie algebras. At the same time, for such algebras there are many 
well-known examples of simple weight modules with infinitely many weights, 
such that all corresponding weight spaces are infinite dimensional 
(see for example \cite{DOF}). In this case 
$\mathrm{soc}_{\mathfrak{n}}(L)$ has infinitely many nonisomorphic
simple submodules, each occurring with infinite multiplicity.

Assume that $(\mathfrak{g},\mathfrak{n})$ is a Whittaker pair and
that there exists a Lie subalgebra $\mathfrak{a}$ of $\mathfrak{g}$
such that $\mathfrak{g}=\mathfrak{a}\oplus \mathfrak{n}$. Then
from the PBW Theorem we have the decomposition 
$U(\mathfrak{g})\cong U(\mathfrak{a})\otimes U(\mathfrak{n})$
(as $U(\mathfrak{a})\text{-}U(\mathfrak{n})$-bimodules) and hence
for any $\lambda\in (\mathfrak{n}/[\mathfrak{n},\mathfrak{n}])^*$
the module $M_{\lambda}$ can be identified with $U(\mathfrak{a})$
as a left $U(\mathfrak{a})$-module via the map
\begin{displaymath}
\begin{array}{ccc}
U(\mathfrak{a}) & \overset{\varphi_{\lambda}}{\longrightarrow} 
& M_{\lambda}\\
u & \mapsto & u\otimes v_{\lambda}.
\end{array}
\end{displaymath}

\begin{proposition}\label{prop35}
Assume that $(\mathfrak{g},\mathfrak{n})$ is a Whittaker pair and
$\mathfrak{g}=\mathfrak{a}\oplus \mathfrak{n}$ for some 
subalgebra $\mathfrak{a}$.  Let 
$\lambda\in (\mathfrak{n}/[\mathfrak{n},\mathfrak{n}])^*$.
Then $\varphi_{\lambda}^{-1}(W_{\lambda}(M_{\lambda}))$ is a subalgebra
of $U(\mathfrak{a})$, isomorphic to 
$\mathrm{End}_{\mathfrak{g}}(M_{\lambda})$.
\end{proposition}

\begin{proof}
We have $\mathrm{End}_{\mathfrak{g}}(M_{\lambda})=
\mathrm{Hom}_{\mathfrak{g}}(M_{\lambda},M_{\lambda})\cong
\mathrm{Hom}_{\mathfrak{n}}(L_{\lambda},M_{\lambda})$ by adjunction
and $\mathrm{Hom}_{\mathfrak{n}}(L_{\lambda},M_{\lambda})\cong
W_{\lambda}(M_{\lambda})$ by the definition of $W_{\lambda}(M_{\lambda})$. 
Now $\varphi_{\lambda}^{-1}$ identifies elements of $W_{\lambda}(M_{\lambda})$
with some elements from $U(\mathfrak{a})$. As $\varphi_{\lambda}$ 
is a homomorphism of $\mathfrak{a}$-modules, this identification is
compatible with the product in $U(\mathfrak{a})$. The claim follows.
\end{proof}

\subsection{Whittaker vectors in completions of highest and lowest
weight modules}
\label{s4.2}

In this subsection we assume that $\mathfrak{g}$ is a Lie algebra 
with a fixed triangular decomposition $\mathfrak{g}=
\mathfrak{n}_-\oplus \mathfrak{h}\oplus \mathfrak{n}_+$ in the sense 
of \cite{MP}. Set $\mathfrak{n}=\mathfrak{n}_+$. Then 
$(\mathfrak{g},\mathfrak{n})$ is a Whittaker pair, see
Example~\ref{ex19}. For $\mu\in \mathfrak{h}^*$ consider
the corresponding Verma module 
\begin{displaymath}
M(\mu):=U(\mathfrak{g}) \bigotimes_{U(\mathfrak{h}\oplus \mathfrak{n}_+)}
\mathbb{C}_{\mu},
\end{displaymath}
where the $U(\mathfrak{h}\oplus \mathfrak{n}_+)$-module structure
on $\mathbb{C}_{\mu}$ is defined as follows: 
$(h+n)(v)=\mu(h)v$, $h\in \mathfrak{h}$, $n\in \mathfrak{n}_+$, 
$v\in \mathbb{C}_{\mu}$, see \cite[Chapter~7]{Di}, \cite[Section~2.3]{MP}.
The module $M(\mu)$ is a highest weight module with highest
weight $\mu$ and we have $M(\mu)=\oplus_{\nu\in \mathfrak{h}^*}M(\mu)_{\nu}$
(see notation of Example~\ref{ex10}). Similarly we can define
the corresponding lowest weight Verma module $N(\mu)$. From the definition
we have that $N(\mu)\cong U(\mathfrak{n})$ as an $\mathfrak{n}$-module
for any $\mu$.

Consider also the completion
\begin{displaymath}
\overline{M(\mu)}:=\prod_{\nu\in \mathfrak{h}^*}M(\mu)_{\nu},
\end{displaymath}
with the induced $\mathfrak{g}$-module structure. Our main result
in this section is the following generalization of \cite[Theorem~3.8]{Ko}.
The proof generalizes the original proof of  \cite[Theorem~3.8]{Ko}.

\begin{theorem}\label{thm45}
Assume that $\mu\in \mathfrak{h}^*$ is such that $M(\mu)$ is simple. 
Then for any $\lambda\in (\mathfrak{n}/[\mathfrak{n},\mathfrak{n}])^*$
we have $\dim W_{\lambda}(\overline{M(\mu)})=1$, in particular, 
$W_{\lambda}(\overline{M(\mu)})$ is a simple $\mathfrak{n}$-module.
\end{theorem}

Note that for Lie algebras with triangular decomposition the assumptions  
of Theorem~\ref{thm45} are generic (see \cite[Chapter~2]{MP}). 

\begin{proof}
Using the canonical antiautomorphism $x\mapsto -x$ of $\mathfrak{g}$
we define on $N(-\mu)^{*}$ the structure of a $\mathfrak{g}$-module.
As $N(-\mu)$ was a lowest weight module with lowest weight
$-\mu$, the weight $\mu$ becomes a highest weight of $N(-\mu)^{*}$.
In particular, from the universal property of Verma modules, we get 
a nonzero $\mathfrak{g}$-module homomorphism, say $\varphi$, from $M(\mu)$ 
to $N(-\mu)^{*}$. This homomorphism is, moreover, injective as the 
module $M(\mu)$ is simple by our assumptions. 

Note that
$\overline{M(\mu)}:=\prod_{\nu\in \mathfrak{h}^*}M(\mu)_{\nu}$
and $N(-\mu)^{*}=\prod_{\nu\in \mathfrak{h}^*}N(-\mu)^*_{-\nu}$.
Recall that $M(\mu)\cong U(\mathfrak{n}_-)$ as an
$\mathfrak{n}_-$-module and $N(-\mu)\cong U(\mathfrak{n}_+)$ 
as an $\mathfrak{n}_+$-module. As $\mathfrak{g}$ has a triangular
decomposition, the corresponding weight spaces in $U(\mathfrak{n}_+)$
and $U(\mathfrak{n}_-)$ have the same dimension. So the 
corresponding weight spaces in $M(\mu)$ and $N(-\mu)$ have the
same dimension as well. From this it
follows that $\varphi$ extends in the obvious way to an isomorphism
$\overline{M(\mu)}\cong N(-\mu)^{*}$. The statement of the theorem
now follows from the following lemma:

\begin{lemma}\label{lem321}
For any $\mu\in\mathfrak{h}^*$ we have 
$\dim W_{\lambda}(N(-\mu)^{*})=1$.
\end{lemma}

\begin{proof}
Consider $N(-\mu)^{*}$ as a $U(\mathfrak{n})$-module, which we 
may identify with $U(\mathfrak{n})^{*}$. 
Denote by $K_{\lambda}$ the kernel of the algebra homomorphism 
$U(\mathfrak{n})\to\mathbb{C}$, induced by $\lambda$. 
For $f\in U(\mathfrak{n})^{*}$,
$x\in \mathfrak{n}$ and $u\in U(\mathfrak{n})$ we have 
$((x-\lambda(x))f)(u)=f(-(x-\lambda(x))u)$. The latter is equal to 
zero for all $u\in U(\mathfrak{n})$ and $x\in \mathfrak{n}$ if 
and only if $f$ annihilates $K_{\lambda}$. However, $K_{\lambda}$
has, by definition, codimension one in $U(\mathfrak{n})$. It follows
that $\dim W_{\lambda}(N(-\mu)^{*})=1$, which completes the proof.
\end{proof}
\end{proof}

Let $\mu\in \mathfrak{h}^*$ 
and $\lambda\in(\mathfrak{n}/[\mathfrak{n},\mathfrak{n}])^*$. 
By  Proposition~\ref{prop30}\eqref{prop30.1}
and Lemma~\ref{lem321}, we have a unique homomorphism
$\varphi:M_{\lambda}\to N(-\mu)^{*}$, which sends $v_{\lambda}$ to a unique
(up to scalar) nonzero vector in $W_{\lambda}(N(-\mu)^{*})$.
Let $L(\lambda,\mu)$ denote the image of this homomorphism.
The following corollary from Theorem~\ref{thm45} and
Lemma~\ref{lem321} gives some evidence for Conjecture~\ref{conj31}:

\begin{corollary}\label{cor43}
Let $\mu\in \mathfrak{h}^*$ and 
$\lambda\in(\mathfrak{n}/[\mathfrak{n},\mathfrak{n}])^*$.
Then $L(\lambda,\mu)$ is a simple $\mathfrak{g}$-module
and satisfies $\dim W_{\lambda}(L(\lambda,\mu))=1$.
\end{corollary}

\begin{proof}
The module $L(\lambda,\mu)$  is generated by $v$ and satisfies 
$\dim W_{\lambda}(L(\lambda,\mu))=1$ by construction
and Lemma~\ref{lem321}. Simplicity of $L(\lambda,\mu)$ now follows 
from Proposition~\ref{prop30}\eqref{prop30.2}. This completes the proof.
\end{proof}

In the case of the Virasoro algebra from the above results one
obtains some new simple modules as we did not have any restriction on 
$\lambda$ compared to \cite{OW}.

\subsection{Annihilators of $L(\lambda,\mu)$}
\label{s4.3}

Here we continue to work in the setup of the previous subsection.
Our main aim here is to prove the following theorem:

\begin{theorem}\label{thm50}
Assume that $\mu\in \mathfrak{h}^*$ is such that $M(\mu)$ is simple
and $\lambda\in(\mathfrak{n}/[\mathfrak{n},\mathfrak{n}])^*$. 
Then $\mathrm{Ann}_{U(\mathfrak{g})}L(\lambda,\mu)=
\mathrm{Ann}_{U(\mathfrak{g})}M(\mu)$.
\end{theorem}

\begin{proof}
Obviously, $\mathrm{Ann}_{U(\mathfrak{g})}M(\mu)\subset 
\mathrm{Ann}_{U(\mathfrak{g})}\overline{M(\mu)}$. As
$L(\lambda,\mu)\subset \overline{M(\mu)}$ by construction, we
thus have $\mathrm{Ann}_{U(\mathfrak{g})}M(\mu)\subset
\mathrm{Ann}_{U(\mathfrak{g})}L(\lambda,\mu)$.

To prove the reversed inclusion we use the arguments from 
the proof of \cite[Theorem~3.9]{Ko}. Let 
$X=\mathrm{Ann}_{U(\mathfrak{g})}M(\mu)$ and
$Y=\mathrm{Ann}_{U(\mathfrak{g})}L(\lambda,\mu)$. Then
both $X$ and $Y$ are stable with respect to the adjoint action of
$\mathfrak{h}$. Assume $u\in Y$ and  $\alpha\in \mathfrak{h}^*$ are 
such that $[h,u]=\alpha(h)u$ for all $h\in \mathfrak{h}$
and that $u\not\in X$. Then there exists some $\xi\in \mathfrak{h}^*$
and $x\in M(\mu)_{\xi}$ such that $ux\neq 0$. Take any  $y\in L(\lambda,\mu)$,
$y\neq 0$, write it as an infinite sum $y=\sum_{\nu}y_{\nu}$ 
of weight vectors, and let $y_{\nu}$ be a nonzero summand. Then 
$y_{\nu}\in M(\mu)$ and hence, by the  simplicity of $M(\mu)$, we 
have $x=ay_{\nu}$ for some weight element $a\in U(\mathfrak{g})$
(of weight $\xi-\nu$). We have $uay=\sum_{\nu}uay_{\nu}=0$
as $ua\in Y$ and $y\in L(\lambda,\mu)$. On the other hand,
the $\alpha+\xi$-component of $uay$ is $ux\neq 0$. The obtained
contradiction shows that $Y\subset X$, which completes the proof.
\end{proof}

\begin{conjecture}\label{conj51}
{\rm  
Let $\mu\in \mathfrak{h}^*$ be such that $M(\mu)$ is simple. Then for 
generic $\lambda\in(\mathfrak{n}/[\mathfrak{n},\mathfrak{n}])^*$ the module
$L(\lambda,\mu)$ is the unique (up to isomorphism) simple
Whittaker module in $\mathfrak{W}_{\mathfrak{n}}^{\mathfrak{g}}(\lambda)$,
whose annihilator coincides with $\mathrm{Ann}_{U(\mathfrak{g})}L(\lambda,\mu)$.
}
\end{conjecture}

Note that, obviously,  $L(\lambda,\mu)\cong L(\lambda',\mu')$ implies 
$\lambda=\lambda'$.

\section{Whittaker modules over the algebra of derivations
of $\mathbb{C}[x_1,x_2,\dots,x_n]$}\label{s5}

\subsection{The algebra $\mathfrak{w}_n$ and its 
decompositions}\label{s5.1}

Denote by $\mathfrak{w}_n$ the Lie algebra of all derivations
of the polynomial ring $\mathbb{C}[x_1,x_2,\dots,x_n]$. 
For $i\in \{1,2,\dots,n\}$ and $\mathbf{m}=(m_1,m_2,\dots,m_n)\in
\{0,1,2,\dots\}^n$ let  
\begin{displaymath}
D_i(\mathbf{m}):=x_1^{m_1}x_2^{m_2}\cdots
x_n^{m_n}\frac{\partial}{\partial x_i}\in \mathfrak{w}_n.
\end{displaymath}
Then $\mathbf{D}:=\{D_i(\mathbf{m})\}$ is a natural basis of 
$\mathfrak{w}_n$. The algebra $\mathfrak{w}_n$ is a simple infinite 
dimensional Lie algebra of Cartan type with the Cartan subalgebra
$\mathfrak{h}$ being the linear span of 
$x_i\frac{\partial}{\partial x_i}$, $i\in \{1,2,\dots,n\}$.

The linear span of the elements 
$x_i\frac{\partial}{\partial x_j}$, $i,j\in \{1,2,\dots,n\}$,
is a Lie subalgebra $\mathfrak{a}$ of $\mathfrak{w}_n$, isomorphic to 
$\mathfrak{gl}_n$. Note that $\mathfrak{h}\subset \mathfrak{a}$ is a
Cartan subalgebra. For the rest we fix some triangular decomposition
\begin{displaymath}
\mathfrak{a}=\mathfrak{n}_-^{\mathfrak{a}}\oplus
\mathfrak{h}\oplus\mathfrak{n}_+^{\mathfrak{a}}.
\end{displaymath}
Denote by $\mathfrak{n}_+$ the linear span of 
$\mathfrak{n}_+^{\mathfrak{a}}$ and the elements 
$\frac{\partial}{\partial x_i}$, $i\in \{1,2,\dots,n\}$.
Denote by $\mathfrak{n}_-$ the linear span of 
$\mathfrak{n}_-^{\mathfrak{a}}$ and all the elements
$D_i(\mathbf{m})$, which are contained in neither $\mathfrak{a}$
nor $\mathfrak{n}_+$. 

\begin{proposition}\label{prop71}
\begin{enumerate}[(i)]
\item\label{prop71.1} We have the decomposition 
$\mathfrak{w}_n=\mathfrak{n}_-\oplus
\mathfrak{h}\oplus\mathfrak{n}_+$ into a direct sum of Lie subalgebras.
\item\label{prop71.2} Both, $(\mathfrak{w}_n,\mathfrak{n}_-)$
and $(\mathfrak{w}_n,\mathfrak{n}_+)$, are Whittaker pairs.
\end{enumerate}
\end{proposition}

\begin{proof}
The claim \eqref{prop71.1} is a straightforward calculation left to 
the reader. To prove the claim \eqref{prop71.2} we have to verify
that both $\mathfrak{n}_-$ and $\mathfrak{n}_+$ are quasi-nilpotent
and that their actions on the corresponding adjoint modules
$\mathfrak{w}_n/\mathfrak{n}_-$ and $\mathfrak{w}_n/\mathfrak{n}_+$
are locally nilpotent. This again can be done by a straightforward
calculation. However, it is easier explained using the geometric picture
of weights.

Let $\{e_i:i=1,\dots,n\}$ be the standard basis of $\mathbb{R}^n$.
Then the weights of the adjoint action of $\mathfrak{h}$
on $\mathfrak{w}_n$ can be viewed as elements from
$\mathbb{R}^n$ in the following way: the element $D_i(\mathbf{m})$
has weight $\mathbf{m}-e_i$. From our
definition of $\mathfrak{n}_+$ and $\mathfrak{n}_-$ it follows
that there exists a hyperplane $H$ in $\mathbb{R}^n$, containing 
$0$ (the weight of $\mathfrak{h}$), such that the weights of
$\mathfrak{n}_+$ and $\mathfrak{n}_-$ belong to different halfspaces with
respect to $H$. Commutation of elements in $\mathfrak{w}_n$ corresponds, as
usual, to the addition  of respective weights  (in fact, the claim \eqref{prop71.1} follows from this observation). As $\mathfrak{n}_+$
is finite dimensional, it follows that it must be nilpotent
(as adding nonzero vectors we eventually would always leave the 
finite set of weights of $\mathfrak{n}_+$ in a finite number of steps).
Similarly one shows that $\mathfrak{n}_-$ is quasi-nilpotent: 
commuting elements from $\mathfrak{n}_-$ we are moving the set of obtained 
weights further and further from the hyperplane $H$. In the limit,
no weights will be left.

Now if we take some point in one of the halfspaces with respect to
$H$ and add to this point vectors from the other halfspace, 
representing weights of $\mathfrak{n}_+$ or $\mathfrak{n}_-$, 
respectively, we would always eventually obtain a point 
from the other halfspace. This shows that the action of $\mathfrak{n}_-$ 
and $\mathfrak{n}_+$  on the respective adjoint modules
$\mathfrak{w}_n/\mathfrak{n}_-$ and $\mathfrak{w}_n/\mathfrak{n}_+$
is locally nilpotent. The claim of the proposition follows.
\end{proof}

We would like to emphasize that
the decomposition given by Proposition~\ref{prop71}\eqref{prop71.1}
is not a triangular decomposition in the sense of \cite{MP} as the
subalgebra $\mathfrak{n}_+$ is finite dimensional while the
subalgebra $\mathfrak{n}_-$ is infinite dimensional.

The algebra $\mathfrak{w}_n$ has a subalgebra $\mathfrak{a}_1$,
spanned by $\mathfrak{a}$, $\mathfrak{n}_+$, and the elements
$\displaystyle x_j\sum_{i=1}^nx_i\frac{\partial}{\partial x_i}$,
$j\in \{1,2,\dots,n\}$.
The algebra $\mathfrak{a}_1$ is isomorphic to $\mathfrak{sl}_{n+1}$
and $\mathfrak{h}\oplus\mathfrak{n}_+$ is a Borel subalgebra of
$\mathfrak{a}_1$.

\subsection{Highest weight and lowest weight 
$\mathfrak{w}_n$-modules}\label{s5.2}

For $\mu\in\mathfrak{h}^*$ consider the simple $\mathfrak{h}$-module
$\mathbb{C}_{\mu}=\mathbb{C}$ with the action given by
$h(v)=\mu(h)v$, $h\in \mathfrak{h}$, $v\in \mathbb{C}_{\mu}$.
Setting $\mathfrak{n}_+\mathbb{C}_{\mu}=0$ we extend 
$\mathbb{C}_{\mu}$ to a $\mathfrak{h}\oplus\mathfrak{n}_+$-module
and can define the corresponding {\em highest weight Verma module}
\begin{displaymath}
M^+(\mu):=U(\mathfrak{w}_n)
\bigotimes_{U(\mathfrak{h}\oplus\mathfrak{n}_+)}\mathbb{C}_{\mu}.
\end{displaymath}
As usual, the module $M^+(\mu)$ has a unique simple quotient,
denoted by $L^+(\mu)$, which is a simple highest weight module
of highest weight $\mu$. As usual, simple highest weight 
modules are classified by their highest weights. Both
$M^+(\mu)$ and $L^+(\mu)$ obviously have finite dimensional
weight spaces.

Setting $\mathfrak{n}_-\mathbb{C}_{\mu}=0$ we extend 
$\mathbb{C}_{\mu}$ to a $\mathfrak{h}\oplus\mathfrak{n}_-$-module
and can define the corresponding {\em lowest weight Verma module}
\begin{displaymath}
M^-(\mu):=U(\mathfrak{w}_n)
\bigotimes_{U(\mathfrak{h}\oplus\mathfrak{n}_-)}\mathbb{C}_{\mu}.
\end{displaymath}
As usual, the module $M^-(\mu)$ has a unique simple quotient,
denoted by $L^-(\mu)$, which is a simple lowest weight module
of lowest weight $\mu$. As usual, simple lowest weight modules 
are classified by their lowest weights. Both
$M^-(\mu)$ and $L^-(\mu)$ obviously have finite dimensional
weight spaces.

Consider the full subcategory $\mathfrak{X}$ of the category of
all $\mathfrak{w}_n$-modules, which consists of all weight
(with respect to $\mathfrak{h}$) modules with finite dimensional
weight spaces. For $M\in \mathfrak{X}$ we can write
$M=\oplus_{\mu\in \mathfrak{h}^*}M_{\mu}$, where each 
$M_{\mu}$ is finite dimensional. Then
$M^{\star}:=\oplus_{\mu\in \mathfrak{h}^*}M_{\mu}^*$,
which is a subspace of the full dual space $M^*$, becomes
a $\mathfrak{w}_n$-module via $(xf)(v):=f(-x(v))$ for
$x\in \mathfrak{w}_n$, $f\in M^{\star}$ and $v\in M$.
As usual, this extends to an exact contravariant and involutive
self-equivalence $\star$ on the category $\mathfrak{X}$.

\begin{proposition}\label{prop81}
For any $\mu\in \mathfrak{h}^*$ we have
\begin{displaymath}
L^+(\mu)^{\star}\cong L^-(-\mu)\text{ and }
L^-(-\mu)^{\star}\cong L^+(\mu).
\end{displaymath}
\end{proposition}

\begin{proof}
As $\star$ is defined using the canonical involution $x\mapsto -x$
on $\mathfrak{w}_n$, it sends a module with highest weight $\mu$
to a module with lowest weight $-\mu$ and vice versa. 
As $\star$ is a self-equivalence, it sends simple modules
to simple modules. The claim follows.
\end{proof}

We emphasize the following immediate corollary, which is 
interesting as the analogous equality is certainly wrong for the 
corresponding Verma modules (since the algebra $\mathfrak{n}_-$ 
is much ``bigger'' than the algebra $\mathfrak{n}_+$):

\begin{corollary}\label{cor879}
For all $\mu,\nu\in \mathfrak{h}^*$ we have 
$\dim L^+(\mu)_{\nu}=\dim L^-(-\mu)_{-\nu}$.
\end{corollary}

\begin{proposition}\label{prop83}
The modules $M^-(\mu)$ are generically irreducible.
\end{proposition}

\begin{proof}
The restriction of $M^-(\mu)$ to $\mathfrak{a}_1$ is a Verma module,
which is generically irreducible.
\end{proof}

Note that the $\mathfrak{w}_n$-module $M^-(\mu)$ may be irreducible even if 
its restriction to $\mathfrak{a}_1$ is reducible. We refer the reader
to \cite{Ru} for more details.

\subsection{Whittaker modules for $\mathfrak{w}_n$}\label{s5.3}

For the Whittaker pair $(\mathfrak{w}_n,\mathfrak{n}_-)$ 
simple Whittaker modules form a subclass of modules,
considered in \cite{Ru}. In this subsection we generalize the 
construction of Whittaker modules from Subsection~\ref{s4.2} to 
obtain simple Whittaker modules for the Whittaker pair
$(\mathfrak{w}_n,\mathfrak{n}_+)$. We set $\mathfrak{n}:=
\mathfrak{n}_+$. Note one big difference with the
setup of  Subsection~\ref{s4.2}: the decomposition of the
algebra $\mathfrak{w}_n$ we work with (Subsection~\ref{s5.1}) 
is not a triangular decomposition in the sense of \cite{MP}.

For $\mu\in \mathfrak{h}^*$ we have a decomposition
$L^+(\mu)=\oplus_{\nu\in \mathfrak{h}^*}L^+(\mu)_{\nu}$.
Consider the corresponding completion
\begin{displaymath}
\overline{L^+(\mu)}=\prod_{\nu\in \mathfrak{h}^*}L^+(\mu)_{\nu} 
\end{displaymath}
of $L^+(\mu)$, which becomes a $\mathfrak{w}_n$-module in the 
natural way.

\begin{theorem}\label{thm91}
Assume that $\mu\in \mathfrak{h}^*$ is such that $M^-(-\mu)$ is simple. 
Then for any $\lambda\in (\mathfrak{n}/[\mathfrak{n},\mathfrak{n}])^*$
we have $\dim W_{\lambda}(\overline{L^+(\mu)})=1$, in particular, 
$W_{\lambda}(\overline{L^+(\mu)})$ is a simple $\mathfrak{n}$-module.
\end{theorem}

\begin{proof}
Mutatis mutandis the proof of Theorem~\ref{thm45}. Note that
we have $M^-(-\mu)=L^-(-\mu)$ under our assumptions, however
$M^+(\mu)\not\cong L^+(\mu)$. Thus in all arguments from the
proof of Theorem~\ref{thm45} one should consider $L^+(\mu)$
instead of the corresponding Verma module.
\end{proof}

For $\mu$ and $\lambda$ as in Theorem~\ref{thm91} let 
$L(\lambda,\mu)$ denote the image of the unique (up to scalar)
nonzero homomorphism from $M_{\lambda}$ to $\overline{L^+(\mu)}$.
Similarly to Corollary~\ref{cor43} we have:

\begin{corollary}\label{cor92}
Assume that $\mu\in \mathfrak{h}^*$ is such that $M^-(-\mu)$ is simple
and $\lambda\in(\mathfrak{n}/[\mathfrak{n},\mathfrak{n}])^*$.
Then $L(\lambda,\mu)$ is a simple $\mathfrak{w}_n$-module
and satisfies $\dim W_{\lambda}(L(\lambda,\mu))=1$.
\end{corollary}

\begin{proof}
Mutatis mutandis the proof of Corollary~\ref{cor43}.
\end{proof}

\begin{corollary}\label{cor93}
Assume that $\mu\in \mathfrak{h}^*$ is such that $M^-(-\mu)$ is simple
and $\lambda\in(\mathfrak{n}/[\mathfrak{n},\mathfrak{n}])^*$.
Then $\mathrm{Ann}_{U(\mathfrak{g})}L(\lambda,\mu)=
\mathrm{Ann}_{U(\mathfrak{g})}L^+(\mu)$.
\end{corollary}

\begin{proof}
Mutatis mutandis the proof of Theorem~\ref{thm50}.
\end{proof}

It is worth to mention  that annihilators of simple highest weight 
$\mathfrak{w}_1$-modules are described in \cite{CM}.
Note that the restriction of the module $L(\lambda,\mu)$, 
constructed in Theorem~\ref{thm91}, to the algebra
$\mathfrak{a}_1$ obviously has a simple socle
(as $\dim W_{\lambda}(L(\lambda,\mu))=1$), which is a simple
Whittaker module for the Whittaker pair 
$(\mathfrak{a}_1,\mathfrak{n}_+)$.
We finish with the following conjecture:

\begin{conjecture}\label{conj95}
{\rm 
The modules $L(\lambda,\mu)$, $\lambda\neq 0$, constructed 
in Theorem~\ref{thm91}, constitute an exhaustive list of simple 
Whittaker modules with $\lambda\neq 0$ for the
Whittaker pair $(\mathfrak{w}_n,\mathfrak{n}_+)$. Moreover,
$L(\lambda,\mu)\cong L(\lambda,\mu')$ if and only if 
$\mathrm{Ann}_{U(\mathfrak{g})}M^-(-\mu)=
\mathrm{Ann}_{U(\mathfrak{g})}M^-(-\mu')$.
}
\end{conjecture}

\subsection{A general existence theorem}\label{s5.4}

The arguments used in the proof of Theorems~\ref{thm45} and \ref{thm91}
can be easily generalized to prove the following general existence theorem:

\begin{theorem}\label{thm911} 
Let $(\mathfrak{g},\mathfrak{n})$ be a Whittaker pair and assume that
$\mathfrak{g}=\mathfrak{a}\oplus\mathfrak{n}$ for some Lie subalgebra
$\mathfrak{a}$. Then for every $\lambda:\mathfrak{n}\to\mathbb{C}$
there exists a (simple) module 
$M\in\mathfrak{W}_{\mathfrak{n}}^{\mathfrak{g}}(\lambda)$ such that 
$W_{\lambda}(M)$ is one-dimensional.
\end{theorem}

\begin{proof}
Let $L$ be any one-dimensional $\mathfrak{a}$-module. Then the 
$\mathfrak{a}$-module $\mathrm{Ind}_{\mathfrak{a}}^{\mathfrak{g}}(L)$
is isomorphic to $\mathfrak{n}$ as $\mathfrak{n}$-module. Analogously
to the proofs of Theorems~\ref{thm45} and \ref{thm91} we get that 
$W_{\lambda}(\mathrm{Ind}_{\mathfrak{a}}^{\mathfrak{g}}(L)^*)$ is
one-dimensional. The claim follows.
\end{proof}

\section{Further examples: Whittaker modules over solvable 
finite dimensional Lie algebras}\label{s6}

\subsection{The two-dimensional solvable algebra}\label{s6.1}

Consider the two-dimensional Lie algebra $\mathfrak{g}$
with the basis $\{a,b\}$ and the Lie bracket given by $[a,b]=b$. 
Let $\mathfrak{n}=\langle b\rangle$, then $(\mathfrak{g},\mathfrak{n})$ 
is a Whittaker pair (Example~\ref{ex1902}). All simple modules over the 
algebra $\mathfrak{g}$ are described in \cite{Bl}, however,
we are also interested in the description of the category of
Whittaker modules (at least those of finite length). 
For $\lambda\in\mathbb{C}$ we consider the standard Whittaker module
$M_{\lambda}=U(\mathfrak{g})/(b-\lambda)$.

\begin{proposition}\label{prop6001}
Let $\lambda\in\mathbb{C}$, $\lambda\neq 0$.
\begin{enumerate}[(i)]
\item\label{prop6001-1} The module
$M_{\lambda}$ is irreducible and hence is a unique (up to isomorphism)
simple object in $\mathfrak{W}_{\mathfrak{n}}^{\mathfrak{g}}(\lambda)$.
\item\label{prop6001-2} We have
\begin{displaymath}
\mathrm{Ext}^i_{\mathfrak{g}}(M_{\lambda},M_{\lambda})\cong
\begin{cases}
\mathbb{C}, &  i=0;\\
0, & i>0,
\end{cases}
\end{displaymath}
in particular, the category 
$\overline{\mathfrak{W}_{\mathfrak{n}}^{\mathfrak{g}}(\lambda)}$
is equivalent to $\mathbb{C}\text{-}\mathrm{mod}$.
\end{enumerate}
\end{proposition}

\begin{proof}
The claim \eqref{prop6001-1} follows from \cite{Bl}. Here is a short argument:
By the PBW Theorem we have $M_{\lambda}\cong \mathbb{C}[a]$ as a 
$\mathbb{C}[a]$-module. In particular, $M_{\lambda}$ is generated by
$1$ as a $\mathbb{C}[a]$-module. The action of $b$ on $\mathbb{C}[a]$ is
given by $b\cdot f(a)=\lambda f(a-1)$, $f(a)\in\mathbb{C}[a]$. Assume
that $M_{\lambda}$ is not simple and let $N$ be a proper submodule
of $M_{\lambda}$. Let $f(a)\in N$ be a nonzero element of minimal
degree. Then $\deg f(a)=k>0$ as $M_{\lambda}$ is generated by scalars
already as a $\mathbb{C}[a]$-module. For the element
$f(a)-\frac{1}{\lambda}b\cdot f(a)\in N$ we have 
\begin{displaymath}
\deg(f(a)-\frac{1}{\lambda}b\cdot f(a))=
\deg(f(a)-f(a-1))=k-1,
\end{displaymath}
a contradiction. The claim \eqref{prop6001-1} follows.

To prove the claim \eqref{prop6001-2} we consider the following
free resolution of $M_{\lambda}$:
\begin{displaymath}
0\to U(\mathfrak{g})\overset{\cdot (b-\lambda)}{\longrightarrow} 
U(\mathfrak{g}) \to U(\mathfrak{g})/(b-\lambda)\cong M_{\lambda}\to 0.
\end{displaymath}
Applying $\mathrm{Hom}_{U(\mathfrak{g})}({}_-,M_{\lambda})$ we get
\begin{equation}\label{eq555}
0\to M_{\lambda}\overset{(b-\lambda)\cdot}{\longrightarrow} 
M_{\lambda}\to 0.
\end{equation}
This implies that $\mathrm{Ext}^i_{\mathfrak{g}}(M_{\lambda},M_{\lambda})=0$
for all $i>1$. Further, we have that
$\mathrm{Ext}^1_{\mathfrak{g}}(M_{\lambda},M_{\lambda})$
is isomorphic to the cokernel of the linear operator 
$(b-\lambda)\cdot $ on $M_{\lambda}$. We claim that this cokernel is
zero. Indeed, for $f(a)\in \mathbb{C}[a]$ we have
\begin{displaymath}
(b-\lambda)\cdot f(a)=\lambda (f(a-1)-f(a)).
\end{displaymath}
As $\lambda\neq 0$, the cokernel of $(b-\lambda)\cdot $ equals the
cokernel of the linear operator $f(a)\mapsto f(a-1)-f(a)$. This cokernel 
is obviously zero. Therefore 
$\mathrm{Ext}^1_{\mathfrak{g}}(M_{\lambda},M_{\lambda})=0$
and hence the category 
$\overline{\mathfrak{W}_{\mathfrak{n}}^{\mathfrak{g}}(\lambda)}$
is semisimple. This proves the claim \eqref{prop6001-2} and 
completes the proof of the proposition.
\end{proof}

For $\mu\in\mathbb{C}$ denote by $L(\mu)$ the one-dimensional 
$\mathfrak{g}$-module given by $b\cdot L(\mu)=0$ and $a\cdot v=\mu v$,
$v\in L(\mu)$. 

\begin{proposition}\label{prop6002}
\begin{enumerate}[(i)]
\item\label{prop6002-1} The modules $\{L(\mu):\mu\in\mathbb{C}\}$
constitute an exhaustive and irredundant  list of pairwise nonisomorphic
simple objects in $\mathfrak{W}_{\mathfrak{n}}^{\mathfrak{g}}(0)$.
\item\label{prop6002-2} For $\mu,\nu\in\mathbb{C}$ we have
\begin{eqnarray}
\mathrm{Ext}^1_{\mathfrak{g}}(L(\mu),L(\nu))&\cong&
\begin{cases}
\mathbb{C}, &  \nu\in\{\mu,\mu+1\};\\
0, & \text{otherwise};
\end{cases}  \label{eq05}\\ 
\mathrm{Ext}^2_{\mathfrak{g}}(L(\mu),L(\nu))&\cong& 
\begin{cases}
\mathbb{C}, &  \nu=\mu+1;\\ 
0, & \text{otherwise};
\end{cases} \label{eq06}\\ 
\mathrm{Ext}^i_{\mathfrak{g}}(L(\mu),L(\nu))&\cong& 0,\quad i>2.
\label{eq07}
\end{eqnarray}
\end{enumerate}
\end{proposition}

\begin{proof}
We have $b\cdot M_0=0$ and hence from Proposition~\ref{prop30}\eqref{prop30.1}
it follows that $b\cdot L=0$ for any simple
$L\in \mathfrak{W}_{\mathfrak{n}}^{\mathfrak{g}}(0)$. Therefore any
simple object in $\mathfrak{W}_{\mathfrak{n}}^{\mathfrak{g}}(0)$
is, in fact, a simple module over the polynomial algebra  $\mathbb{C}[a]$.
The claim \eqref{prop6002-1} follows.

For $\mu\in\mathbb{C}$ it is easy to see that the following is a 
free resolution of the module $L(\mu)$:
\begin{multline*}
0\to U(\mathfrak{g})\overset{\left(\begin{array}{c}
\text{\tiny $\cdot(a$-$\mu$-$1)$}\\
\text{\tiny $\cdot b$}\end{array}\right)}
{\longrightarrow} U(\mathfrak{g})\oplus
U(\mathfrak{g})\overset{(-\cdot b,\cdot (a-\mu))}{\longrightarrow} 
U(\mathfrak{g}) \to\\\to U(\mathfrak{g})/(a-\mu,b)\cong L(\mu)\to 0.
\end{multline*}
In particular, all extensions of degree three and higher vanish
and the formula \eqref{eq07} follows.
Applying $\mathrm{Hom}_{U(\mathfrak{g})}({}_-,L(\nu))$, 
$\nu\in\mathbb{C}$, we get
\begin{equation}\label{eq9}
0\to L(\nu)\overset{\left(\begin{array}{c}
\text{\tiny $-b\cdot $}\\
\text{\tiny $(a$-$\mu)\cdot$}\end{array}\right)}
{\longrightarrow} L(\nu)\oplus L(\nu)
\overset{((a-\mu-1)\cdot ,b\cdot)}{\longrightarrow} L(\nu)\to 0.
\end{equation}
As $L(\nu)\cong \mathbb{C}$ and $b\cdot L(\nu)=0$, in the case 
$\nu\not\in\{\mu,\mu+1\}$ we immediately obtain that \eqref{eq9}
is exact. If $\nu=\mu$, we have one dimensional homologies in degrees
zero and one.  
If $\nu=\mu+1$, we have one dimensional homologies in degrees
one  and two. This gives the formulae \eqref{eq05} and \eqref{eq06}.
The claim  \eqref{prop6002-2} follows and the proof is complete.
\end{proof}

Now we would like to determine a decomposition of the category
$\overline{\mathfrak{W}_{\mathfrak{n}}^{\mathfrak{g}}(0)}$
into indecomposable subcategories. For $\xi\in\mathbb{C}/\mathbb{Z}$ 
denote by  $\overline{\mathfrak{W}_{\mathfrak{n}}^{\mathfrak{g}}(0)}_{\xi}$
the full subcategory of 
$\overline{\mathfrak{W}_{\mathfrak{n}}^{\mathfrak{g}}(0)}$,
which consists of all modules, whose simple subquotients have the form
$L(\mu)$, $\mu\in\xi$.

\begin{theorem}\label{thm6003}
\begin{enumerate}[(i)]
\item\label{thm6003.1} $\displaystyle
\overline{\mathfrak{W}_{\mathfrak{n}}^{\mathfrak{g}}(0)}\cong
\bigoplus_{\xi\in\mathbb{C}/\mathbb{Z}}
\overline{\mathfrak{W}_{\mathfrak{n}}^{\mathfrak{g}}(0)}_{\xi}$.
\item\label{thm6003.2} Each 
$\overline{\mathfrak{W}_{\mathfrak{n}}^{\mathfrak{g}}(0)}_{\xi}$
is equivalent to the category of finite-dimensional modules over
the following quiver:
\begin{equation}\label{eq12}
\xymatrix{ 
\dots\ar[r]_{\mathtt{b}_{i-1}}& 
\bullet\ar[r]_{\mathtt{b}_i}\ar@(ul,ur)[]^{\mathtt{a}_i} & 
\bullet\ar[r]_{\mathtt{b}_{i+1}}\ar@(ul,ur)[]^{\mathtt{a}_{i+1}} 
& \bullet\ar[r]_{\mathtt{b}_{i+2}}\ar@(ul,ur)[]^{\mathtt{a}_{i+2}} &\dots
},
\end{equation}
with relations $\mathtt{b}_i\mathtt{a}_i=\mathtt{a}_{i+1}\mathtt{b}_i$
for all $i$, where every $\mathtt{a}_i$ acts locally nilpotent.
\end{enumerate}
\end{theorem}

\begin{proof}
If $\xi\in\mathbb{C}/\mathbb{Z}$ and $\mu\not\in\xi$, then from 
Proposition~\ref{prop6002}\eqref{prop6002-2} we have that the first
extension from $L(\mu)$ to any $L(\nu)$, $\nu\in \xi$, an vice versa,
vanishes. The claim  \eqref{thm6003.1} follows.

Fix now $\xi\in\mathbb{C}/\mathbb{Z}$. As $\mathfrak{n}$ is finite
dimensional, the category 
$\mathfrak{W}_{\mathfrak{n}}^{\mathfrak{g}}$
is extension closed in $\mathfrak{g}\text{-}\mathrm{Mod}$ by
Proposition~\ref{prop1}. Hence, by 
Proposition~\ref{prop6002}\eqref{prop6002-2},
the quiver given in \eqref{thm6003.2} is the ext-quiver of
$\overline{\mathfrak{W}_{\mathfrak{n}}^{\mathfrak{g}}(0)}_{\xi}$.
The relation 
\begin{equation}\label{eq912}
\mathtt{b}_i\mathtt{a}_i=\mathtt{a}_{i+1}\mathtt{b}_i
\end{equation}
follows immediately from the relation $ab=b(a+1)$ in $U(\mathfrak{g})$.

Assume that $w=0$ is a new relation for
$\overline{\mathfrak{W}_{\mathfrak{n}}^{\mathfrak{g}}(0)}_{\xi}$, 
which does not follow from the relations \eqref{eq912}, and that the 
maximal degree of a monomial in $w$ is $k-1$ for some $k\in\mathbb{N}$.

Consider the full subcategory $\mathfrak{X}$ of
$\overline{\mathfrak{W}_{\mathfrak{n}}^{\mathfrak{g}}(0)}_{\xi}$,
consisting of all $M$ such that $b^k\cdot M=0$ and such that 
the minimal polynomial of the action of $a$ on $M$ has roots of
multiplicities at most $k$.

For $k\in\mathbb{Z}$ and $\mu\in \xi$ consider the $\mathfrak{g}$-module
\begin{displaymath}
V_k(\mu)=U(\mathfrak{g})/((a-\mu)^k,b^k)\in
\overline{\mathfrak{W}_{\mathfrak{n}}^{\mathfrak{g}}(0)}_{\xi}.
\end{displaymath}
It is easy to see that the endomorphism algebra of this module
is isomorphic to $\mathbb{C}[x]/(x^k)$ (acting via multiplication with
$a-\mu$ from the right). In particular, the module
$V_k(\mu)$ is indecomposable for all $k$ and is generated by $1$. 
It is easy to see that $V_k(\mu)$ is, in fact, an
indecomposable projective in $\mathfrak{X}$.
From the PBW Theorem we have that $\dim V_k(\mu)=k^2$. 

Consider the quotient $A$ of our quiver algebra 
\eqref{eq12}, given by \eqref{eq912} and additional relations
(for all $i$) 
\begin{equation}\label{eq914}
\mathtt{a}^k_i=0,\quad  
\mathtt{b}_{i+k-1}\cdots\mathtt{b}_{i+1}\mathtt{b}_i=0.
\end{equation}
Then all indecomposable projective $A$-modules have dimension
$k^2$. If we would add the additional relation $w=0$ (which is not
a consequence of \eqref{eq912} and \eqref{eq914} by our
choice of $k$), the dimension of indecomposable projectives
$A$-modules would decrease. This, however, contradicts
the previous paragraph.  The claim of the theorem follows.
\end{proof}

The quiver algebra, described in Theorem~\ref{thm6003}\eqref{thm6003.2}, 
is Koszul (we refer the reader to \cite{MOS} for Koszul theory for
algebras with infinitely many simple modules). In particular, there is a 
graded version of the category  
$\overline{\mathfrak{W}_{\mathfrak{n}}^{\mathfrak{g}}(0)}_{\xi}$, which
is equivalent to the category of finite dimensional modules over
the following quiver subject to the relations that all squares commutes:
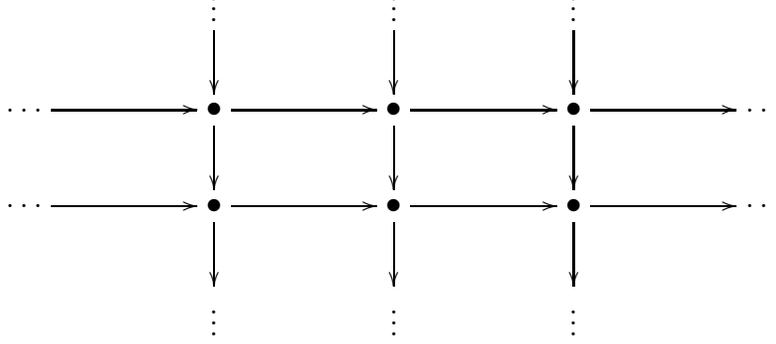
\begin{displaymath}
\xymatrix{ 
 &&\vdots\ar[d] &&\vdots\ar[d] && \vdots\ar[d]&& \\
\dots\ar[rr] &&\bullet\ar[d]\ar[rr] &&\bullet\ar[d]\ar[rr]  
&& \bullet\ar[d]\ar[rr] && \dots\\
\dots \ar[rr]&&\bullet\ar[d]\ar[rr] &&\bullet\ar[d]\ar[rr]  
&& \bullet\ar[d]\ar[rr] && \dots\\
 &&\vdots &&\vdots && \vdots&& 
}
\end{displaymath}
This is interesting as a priori we do not have any grading on 
$\overline{\mathfrak{W}_{\mathfrak{n}}^{\mathfrak{g}}(0)}_{\xi}$.

\subsection{A three-dimensional nilpotent algebra}\label{s6.2}

Let $\mathfrak{g}$ denote the three-dimensional Lie algebra with the
basis $\{a,b,c\}$, where the Lie bracket is given by
\begin{displaymath}
[a,b]=c,\quad [a,c]=0,\quad [b,c]=0.
\end{displaymath}
Let $\mathfrak{n}=\langle c\rangle$, then 
$(\mathfrak{g},\mathfrak{n})$ is a Whittaker pair (Example~\ref{ex1902}).
As $c$ is central it acts as a scalar on any simple 
$\mathfrak{g}$-module (\cite[2.6.8]{Di}). Hence any simple
$\mathfrak{g}$-module is a Whittaker module for the Whittaker
pair $(\mathfrak{g},\mathfrak{n})$ (in particular, it might
have many Whittaker vectors, compare with Conjecture~\ref{conj31}). 
For every $\theta\in\mathbb{C}$
the quotient $U(\mathfrak{g})/(c-\theta)$ is isomorphic to the
first Weyl algebra $\mathrm{A}_1$. Simple modules over this algebra 
are described in \cite{Bl} and the blocks for this algebra
seem to be very complicated, see for example \cite{Bav}. 

\subsection{Borel subalgebras}\label{s6.3}

Let $\mathfrak{g}$ be a semi-simple finite dimensional Lie algebra
with a fixed triangular decomposition $\mathfrak{g}=
\mathfrak{n}_-\oplus \mathfrak{h}\oplus\mathfrak{n}_+$. Then
$\mathfrak{b}=\mathfrak{h}\oplus\mathfrak{n}_+$ is a Borel subalgebra
of $\mathfrak{g}$ and $(\mathfrak{b},\mathfrak{n})$,
where $\mathfrak{n}=\mathfrak{n}_+$, is a 
Whittaker pair (Example~\ref{ex1902}). Note that the Whittaker pair 
considered in Subsection~\ref{s6.1} was a special case of the
present situation ($\mathfrak{g}=\mathfrak{sl}_2$). 
Let $\Delta_+$ be the set of positive
roots for $\mathfrak{g}$ (that is roots for $\mathfrak{n}$) and $\pi$ be 
the set of simple roots. For $\alpha\in\Delta_+$ we fix some nonzero
root vector $e_{\alpha}\in \mathfrak{g}_{\alpha}\subset \mathfrak{n}$. 
Let $h_{\alpha}$, $\alpha\in\pi$, be the basis of $\mathfrak{h}$, 
dual to  $\pi$. We have
$[\mathfrak{n},\mathfrak{n}]=\langle e_{\alpha}:\alpha\in 
\Delta_+\setminus \pi\rangle$. In particular, every Lie algebra
homomorphism $\lambda:\mathfrak{n}
\to\mathbb{C}$ satisfies $\lambda(e_{\alpha})=0$,
$\alpha\in \Delta_+\setminus \pi$, and is uniquely defined by 
$\lambda_{\alpha}=\lambda(e_{\alpha})$, $\alpha\in \pi$.

For $\lambda:\mathfrak{n}\to\mathbb{C}$ set
$\pi_{\lambda}=\{\alpha\in\pi:\lambda_{\alpha}\neq 0\}$. Define
\begin{displaymath}
\mathfrak{h}^{\lambda}=\langle h_{\alpha}:\alpha\in\pi_{\lambda}\rangle,
\quad\quad
\mathfrak{h}_{\lambda}=\langle h_{\alpha}:\alpha\in\pi\setminus
\pi_{\lambda}\rangle.
\end{displaymath}
For $\lambda:\mathfrak{n}\to\mathbb{C}$ and $\mu\in\mathfrak{h}_{\lambda}^*$
let $I_{\lambda,\mu}$ denote the left ideal of $U(\mathfrak{b})$ generated
by the elements $e_{\alpha}-\lambda_{\alpha}$, $\alpha\in\pi$,
and $h-\mu(h)$, $h\in \mathfrak{h}_{\lambda}$. Set $L_{\lambda,\mu}=
U(\mathfrak{b})/I_{\lambda,\mu}$. From the PBW Theorem we have 
that $U(\mathfrak{b})/I_{\lambda,\mu} \cong\mathbb{C}
[\mathfrak{h}^{\lambda}]$ as a $\mathbb{C}[\mathfrak{h}^{\lambda}]$-module.

\begin{proposition}\label{prop701}
The modules $\{L_{\lambda,\mu}:\lambda:\mathfrak{n}\to\mathbb{C},
\mu\in\mathfrak{h}_{\lambda}^*\}$ constitute an exhaustive and irredundant
list of simple modules in $\mathfrak{W}_{\mathfrak{n}}^{\mathfrak{b}}$.
\end{proposition}

\begin{proof}
From the definition it follows easily that the canonical generator of
$L_{\lambda,\mu}$ is a unique (up to scalar) Whittaker element of
$L_{\lambda,\mu}$. Hence simplicity of $L_{\lambda,\mu}$ follows
from Proposition~\ref{prop30}\eqref{prop30.2}. That these modules
are pairwise non-isomorphic follows directly from the definition.

For a fixed $\lambda:\mathfrak{n}\to\mathbb{C}$ let 
$L\in \mathfrak{W}_{\mathfrak{n}}^{\mathfrak{b}}(\lambda)$ be simple.
Then $L$ is a quotient of $M_{\lambda}$ by
Proposition~\ref{prop30}\eqref{prop30.1}, in particular,
$e_{\alpha}L=0$ for any $\alpha\in \Delta_{+}\setminus\pi_{\lambda}$.
Hence that action of $U(\mathfrak{h}_{\lambda})$ on $L$ 
commutes with the action of the whole $U(\mathfrak{b})$ and thus 
gives endomorphisms of $L$. However, every endomorphism of a
simple module reduces to scalars by \cite[2.6.4]{Di}, so every
element of $U(\mathfrak{h}_{\lambda})$ acts on $L$ as a scalar,
say the one, given by $\mu\in\mathfrak{h}_{\lambda}^*$.
This implies that $I_{\lambda,\mu}$ annihilates the Whittaker vector of
$L$. The latter means that $L$ is a nonzero quotient of $L_{\lambda,\mu}$
and hence $L\cong L_{\lambda,\mu}$ as $L_{\lambda,\mu}$ is simple.
This completes the proof.
\end{proof}

For every $\alpha\in \pi$ the space $\mathfrak{a}_{\alpha}=\langle 
h_{\alpha},e_{\alpha}\rangle$ is a Lie subalgebra of $\mathfrak{b}$,
isomorphic to the Lie algebra $\mathfrak{g}$ from 
Subsection~\ref{s6.1} and we have
\begin{displaymath}
\bigoplus_{\alpha\in \pi}\mathfrak{a}_{\alpha}\cong\mathfrak{a}:=
\mathfrak{b}/[\mathfrak{n},\mathfrak{n}].
\end{displaymath}
For $\lambda:\mathfrak{n}\to\mathbb{C}$ and $\mu\in\mathfrak{h}_{\lambda}^*$
we have $[\mathfrak{n},\mathfrak{n}]M_{\lambda}=0$ from the definition
and hence $[\mathfrak{n},\mathfrak{n}]L_{\lambda,\mu}=0$ as well.
This makes $L_{\lambda,\mu}$ an $\mathfrak{a}$-module and we have
an obvious isomorphism of $\mathfrak{a}$-modules:
\begin{equation}\label{eq198}
L_{\lambda,\mu}\cong\bigotimes_{\alpha\in \pi_{\lambda}}
L^{\mathfrak{a}_{\alpha}}_{\lambda_{\alpha}}\otimes 
\bigotimes_{\alpha\in \pi\setminus\pi_{\lambda}}
L^{\mathfrak{a}_{\alpha}}(\mu(h_{\alpha}))
\end{equation}
(here the superscript $\mathfrak{a}_{\alpha}$ means that the module in
question is an $\mathfrak{a}_{\alpha}$-module and the notation is as
in  Subsection~\ref{s6.1}). This isomorphism extends to 
$\mathfrak{b}$-modules using the trivial  action of $[\mathfrak{n},\mathfrak{n}]$. To simplify notation we set
\begin{displaymath}
L^{\mathfrak{a}_{\alpha}}_{\lambda,\mu}:=
\begin{cases}
L^{\mathfrak{a}_{\alpha}}_{\lambda_{\alpha}}, & \alpha\in \pi_{\lambda};\\
L^{\mathfrak{a}_{\alpha}}(\mu(h_{\alpha})), & 
\alpha\in \pi\setminus\pi_{\lambda}.
\end{cases}
\end{displaymath}

\begin{proposition}\label{prop7012}
Let $\lambda:\mathfrak{n}\to\mathbb{C}$,
$\mu\in\mathfrak{h}_{\lambda}^*$, $\lambda':\mathfrak{n}\to\mathbb{C}$,
$\mu'\in\mathfrak{h}_{\lambda}^*$. Then
\begin{displaymath}
\mathrm{Ext}^k_{\mathfrak{a}}(L_{\lambda,\mu},L_{\lambda',\mu'})\cong
\sum_{\text{\tiny$\begin{array}{c}(i_{\alpha})\in\{0,1,2,\dots\}^{\pi}\\
\sum_{\alpha} i_{\alpha}=k\end{array}$}}
\bigotimes_{\alpha\in\pi}
\mathrm{Ext}^{i_{\alpha}}_{\mathfrak{a}_{\alpha}}
(L^{\mathfrak{a}_{\alpha}}_{\lambda,\mu},
L^{\mathfrak{a}_{\alpha}}_{\lambda',\mu'}),
\end{displaymath}
in particular, all these extension spaces are finite dimensional.
\end{proposition}

We note that all extension spaces on the right hand side of the 
above formula are explicitly described in Subsection~\ref{s6.1}.

\begin{proof}
As $\mathfrak{a}\cong\oplus_{\alpha\in\pi}\mathfrak{a}_{\alpha}$
and $U(\mathfrak{a})\cong\otimes_{\alpha\in\pi}U(\mathfrak{a}_{\alpha})$,
the claim follows from \eqref{eq198} and the K{\"u}nneth formula.
\end{proof}

\begin{corollary}\label{cor7015}
All vector spaces 
$\mathrm{Ext}^1_{\mathfrak{b}}(L_{\lambda,\mu},L_{\lambda',\mu'})$
are finite dimensional.
\end{corollary}

\begin{proof}
To compute the extension spaces in question we use the classical cohomology
of Lie algebras (\cite[Chapter~XIII]{CE}). Consider the Lie algebra 
$\mathfrak{c}=\mathfrak{n}\oplus\mathfrak{h}_{\lambda}$ and observe that
$L_{\lambda,\mu}$ is induced from the one-dimensional 
$\mathfrak{c}$-module $\mathtt{L}_{\lambda,\mu}$, given by 
$\lambda$ and $\mu$, by construction. We denote by 
$\nu:\mathfrak{c}\to \mathbb{C}$ the Lie algebra homomorphism, defining
$\mathtt{L}_{\lambda,\mu}$. Consider the first three steps
of the free resolution of the trivial $\mathfrak{c}$-module $\mathbb{C}$:
\begin{equation}\label{eq101}
U(\mathfrak{c})\otimes\mathfrak{c}\wedge  \mathfrak{c}
\overset{\gamma_1}{\rightarrow} U(\mathfrak{c})\otimes\mathfrak{c} 
\overset{\beta_1}{\rightarrow} U(\mathfrak{c})\otimes\mathbb{C} 
\overset{\alpha_1}{\rightarrow}\mathbb{C}\to 0,
\end{equation}
where the map $\alpha$ is given by $1\otimes 1\mapsto 1$;
the map $\beta$ is given by $1\otimes x\mapsto x$, $x\in \mathfrak{c}$;
and the map $\gamma$ is given by 
\begin{displaymath}
1\otimes x\wedge y\mapsto x\otimes y - y\otimes x -1\otimes [x,y], 
\quad\quad x,y\in \mathfrak{c}.
\end{displaymath}
Tensoring \eqref{eq101} with the $\mathfrak{c}$-module
$\mathtt{L}_{\lambda,\mu}$ (over $\mathbb{C}$) and further
with $U(\mathfrak{b})$ over $U(\mathfrak{c})$
(which is exact by the PBW Theorem), we obtain a free 
$U(\mathfrak{b})$-resolution of $L_{\lambda,\mu}\cong
U(\mathfrak{b})\otimes_{U(\mathfrak{c})}
\left(\mathbb{C}\otimes\mathtt{L}_{\lambda,\mu}\right)$ 
as follows:
\begin{multline}\label{eq103}
U(\mathfrak{b})\otimes_{U(\mathfrak{c})}\left(U(\mathfrak{c})
\otimes\mathfrak{c}\wedge  \mathfrak{c}
\otimes\mathtt{L}_{\lambda,\mu}\right)
\overset{\gamma_2}{\rightarrow} U(\mathfrak{b})\otimes_{U(\mathfrak{c})}
\left(U(\mathfrak{c})
\otimes\mathfrak{c}\otimes\mathtt{L}_{\lambda,\mu}\right) 
\overset{\beta_2}{\rightarrow}\\
\overset{\beta_2}{\rightarrow} 
U(\mathfrak{b})\otimes_{U(\mathfrak{c})}\left(U(\mathfrak{c})
\otimes\mathbb{C}
\otimes\mathtt{L}_{\lambda,\mu}\right) 
\overset{\alpha_2}{\rightarrow}U(\mathfrak{b})\otimes_{U(\mathfrak{c})}
\left(\mathbb{C}\otimes\mathtt{L}_{\lambda,\mu}\right)\to 0,
\end{multline}
where $x_2=\mathrm{id}\otimes (x_1\otimes\mathrm{id})$,
$x\in\{\alpha,\beta,\gamma\}$.

Now we would like to  apply 
$\mathrm{Hom}_{\mathfrak{b}}({}_-,L_{\lambda',\mu'})$ to 
\eqref{eq103} (omitting the term 
$U(\mathfrak{b})\otimes_{U(\mathfrak{c})}
\left(\mathbb{C}\otimes\mathtt{L}_{\lambda,\mu}\right)$).
As components of \eqref{eq103} are free $\mathfrak{c}$-modules of finite
rank, the result will be a complex, every component of which is a direct sum 
of some copies of $L_{\lambda',\mu'}$. To be able to write down the maps
explicitly (which will be necessary for our computations), 
we would need some notation and a rewritten version  of \eqref{eq103}.

Choose some ordered basis $\mathbf{b}=\{b_1,\dots,b_k\}$ 
of $\mathfrak{c}$ consisting of the $e_{\alpha}$'s for  $\alpha\in\Delta_+$, 
and the $h_{\alpha}$'s for $\alpha\in\pi\setminus\pi_{\lambda}$. We 
can identify 
\begin{displaymath}
\begin{array}{ccc}
U(\mathfrak{b})\otimes_{U(\mathfrak{c})}
\left(\mathbb{C}\otimes\mathtt{L}_{\lambda,\mu}\right)
&\text{ with }& U(\mathfrak{b})/I_{\lambda,\mu},\\
U(\mathfrak{b})\otimes_{U(\mathfrak{c})}\left(U(\mathfrak{c})
\otimes\mathbb{C}
\otimes\mathtt{L}_{\lambda,\mu}\right)
&\text{ with }& U(\mathfrak{b}),\\
U(\mathfrak{b})\otimes_{U(\mathfrak{c})}
\left(U(\mathfrak{c})
\otimes\mathfrak{c}\otimes\mathtt{L}_{\lambda,\mu}\right) 
&\text{ with }& \displaystyle \bigoplus_{i=1,\dots,k} U(\mathfrak{b}),\\
U(\mathfrak{b})\otimes_{U(\mathfrak{c})}\left(U(\mathfrak{c})
\otimes\mathfrak{c}\wedge  \mathfrak{c}
\otimes\mathtt{L}_{\lambda,\mu}\right)
&\text{ with }& \displaystyle \bigoplus_{1\leq i<j\leq k} U(\mathfrak{b}),\\
\end{array}
\end{displaymath}
such that \eqref{eq103} becomes
\begin{equation}\label{eq104}
\bigoplus_{1\leq i<j\leq k} U(\mathfrak{b})
\overset{\gamma_3}{\rightarrow}
\bigoplus_{i=1,\dots,k} U(\mathfrak{b})
\overset{\beta_3}{\rightarrow}
U(\mathfrak{b}) 
\overset{\alpha_3}{\rightarrow}
U(\mathfrak{b})/I_{\lambda,\mu}\to 0,
\end{equation}
where $\alpha_3$ is the natural projection; $\beta_3$ is given by
\begin{displaymath}
B=(\cdot(b_1-\nu(b_1)),\dots,\cdot(b_i-\nu(b_i)))
\end{displaymath}
(the starting dot means ``right multiplication''), and $\gamma_3$ is given by
the matrix $C=(\cdot c_{ij,s})_{1\leq i<j\leq k}^{s=1,\dots,k}$ such that
\begin{equation}\label{eq106}
c_{ij,s}=
\begin{cases}
b_j-\nu(b_j), & s=i, [b_i,b_j]\neq cb_i\text{ for any }
c\in\mathbb{C}\setminus\{0\};\\
-b_i+\nu(b_i), & s=j, [b_i,b_j]\neq cb_j\text{ for any }
c\in\mathbb{C}\setminus\{0\};\\
b_j-\nu(b_j)+c, & s=i, [b_i,b_j]= cb_i\text{ for some }
c\in\mathbb{C}\setminus\{0\};\\
-b_i+\nu(b_i)-c, & s=j, [b_i,b_j]= cb_j\text{ for some }
c\in\mathbb{C}\setminus\{0\};\\
\frac{1}{c}, & b_s=c[b_i,b_j], c\in \mathbb{C}\setminus\{0\},
s\neq i,j;\\
0, & \text{otherwise}. 
\end{cases}
\end{equation}

Applying $\mathrm{Hom}_{\mathfrak{b}}({}_-,L_{\lambda',\mu'})$ to 
\eqref{eq104} (omitting the non-free term $U(\mathfrak{b})/I_{\lambda,\mu}$), 
we obtain the following complex:
\begin{equation}\label{eq105}
0\to  L_{\lambda',\mu'} \overset{\beta_4}{\rightarrow} 
\bigoplus_{i=1,\dots,k} L_{\lambda',\mu'} 
\overset{\gamma_4}{\rightarrow}
\bigoplus_{1\leq i<j\leq k} L_{\lambda',\mu'},
\end{equation}
where the maps $\beta_4$ and $\gamma_4$ are given by the matrices
$B^T$ and $C^T$, respectively (and the left instead of the right
multiplication with the elements of these matrices).

Now we have to estimate the dimension of the first homology of the
complex \eqref{eq105}. Split the direct sum 
$\bigoplus_{i=1,\dots,k} L_{\lambda',\mu'} $ into two parts: the
first one, $X$, corresponding to $b_i\in[\mathfrak{n},\mathfrak{n}]$,
and the rest, $Y$. We will need the following lemma:

\begin{lemma}\label{lem135}
For all $m\geq 0$ the vectorspace 
$\displaystyle \bigcap_{\alpha\in\pi_{\lambda'}}
\mathrm{Ker}(e_{\alpha}-\lambda'(\alpha))^m$ of $L_{\lambda',\mu'}$
is finite dimensional.
\end{lemma}

\begin{proof}
In the case  $|\pi_{\lambda'}|=1$ this follows from the description of the 
action of $e_{\alpha}-\lambda'(\alpha)$ in Subsection~\ref{s6.1}.
As the general case is a tensor power of the case considered in 
Subsection~\ref{s6.1} (see \eqref{eq198}), the general claim follows 
from the fact that 
the tensor product of finite dimensional spaces is finite dimensional.
\end{proof}

As $[\mathfrak{n},\mathfrak{n}]L_{\lambda',\mu'}=0$, the image of 
$\beta_4$ belongs to $Y$. Let now $x=(x_i)
\in \bigoplus_{i=1,\dots,k} L_{\lambda',\mu'}$ and assume that
$x\not\in Y$. Then $x_i\neq 0$ for some component $i$ from $X$. 
Assume for the moment that $b_i$ is central in 
$\mathfrak{n}$. Let $b_j\in\mathbf{b}$, $j\neq i$. 
Then from \eqref{eq106} one obtains that the $ij$-th (or the $ji$-th, 
depending on the ordering of $i$ and $j$) component of 
$\gamma_4(x)$ equals either $\pm(b_j-\nu(b_j))x_i$ or 
$\pm(b_j-\nu(b_j)+c)x_i$. 
Hence $\gamma_4(x)=0$ implies that $x_i$ must be a Whittaker 
vector in $L_{\lambda',\mu'}$ (note that the subspace of all Whittaker
vectors in $L_{\lambda',\mu'}$ is one-dimensional). If $b_i$ is such that 
$[b_i,\mathfrak{n}]$ is in the center of 
$\mathfrak{n}$, then  we can apply similar arguments 
and get (from \eqref{eq106}) that either $\pm(b_j-\nu(b_j))x_i=0$  
(resp. $\pm(b_j-\nu(b_j)+c)x_i=0$) or $\pm(b_j-\nu(b_j))x_i$ 
(resp. $\pm(b_j-\nu(b_j)+c)x_i$) is proportional to the component of
$[b_i,b_j]$, in which case both $b_i,b_j\in \mathfrak{n}$. The latter
means that $[b_i,b_j]$ is central in $\mathfrak{n}$ and thus 
$\pm(b_j-\nu(b_j))x_i$ (resp. $\pm(b_j-\nu(b_j)+c)x_i$) is proportional
to the Whittaker vector in $L_{\lambda',\mu'}$. This means that
$x_i$ must be in the subspace 
$\displaystyle \bigcap_{\alpha\in\pi_{\lambda'}}
\mathrm{Ker}(e_{\alpha}-\lambda'(\alpha))^2$ of $L_{\lambda',\mu'}$,
which is finite dimensional by Lemma~\ref{lem135}. Proceeding by induction,
we obtain that every $x_i$ for a component from $X$ belongs to some 
fixed subspace of the form
$\displaystyle \bigcap_{\alpha\in\pi_{\lambda'}}
\mathrm{Ker}(e_{\alpha}-\lambda'(\alpha))^m$, which is 
finite  dimensional by Lemma~\ref{lem135}.

The above means that $\mathrm{Ker}(\gamma_4)/
(\mathrm{Ker}(\gamma_4)\cap Y)$
is finite dimensional. At the same time, the first homology of 
\eqref{eq105} coming from $Y$ corresponds exactly to 
$\mathrm{Ext}^1_{\mathfrak{a}}(L_{\lambda,\mu},L_{\lambda',\mu'})$
by the same construction as in the first part of the proof, applied 
to the algebra $\mathfrak{a}$. This part is finite dimensional by 
Proposition~\ref{prop7012}. It follows that the whole first homology of 
\eqref{eq105}  is finite dimensional, which completes the proof.
\end{proof}

The natural inclusion $\mathfrak{h}_{\lambda}\hookrightarrow\mathfrak{h}$
induces the natural projection $\mathfrak{h}^*\tto 
\mathfrak{h}_{\lambda}^*$. Let $G$ denote the image of the abelian group
$\mathbb{Z}\Delta$ under this projection. For $\xi\in 
\mathfrak{h}_{\lambda}^*/G$ denote by 
$\overline{\mathfrak{W}_{\mathfrak{n}}^{\mathfrak{b}}(\lambda)}_{\xi}$
the Serre subcategory of 
$\overline{\mathfrak{W}_{\mathfrak{n}}^{\mathfrak{b}}(\lambda)}$,
generated by $L_{\lambda,\mu}$, $\mu\in \xi$. Using the usual theory of 
weight and generalized weight modules one easily proves the 
following block decomposition for the category 
$\overline{\mathfrak{W}_{\mathfrak{n}}^{\mathfrak{b}}(\lambda)}$:
\begin{displaymath}
\overline{\mathfrak{W}_{\mathfrak{n}}^{\mathfrak{b}}(\lambda)}=
\bigoplus_{\xi\in 
\mathfrak{h}_{\lambda}^*/G} 
\overline{\mathfrak{W}_{\mathfrak{n}}^{\mathfrak{b}}(\lambda)}_{\xi}.
\end{displaymath}
An interesting question is to describe 
$\overline{\mathfrak{W}_{\mathfrak{n}}^{\mathfrak{b}}(\lambda)}_{\xi}$
via quiver and relations similarly to Theorem~\ref{thm6003}\eqref{thm6003.2}.
Note that by Proposition~\ref{prop1} we can use
Corollary~\ref{cor7015} to compute extensions in
$\overline{\mathfrak{W}_{\mathfrak{n}}^{\mathfrak{b}}(\lambda)}_{\xi}$. 
By this Corollary, the ext-quiver of 
$\overline{\mathfrak{W}_{\mathfrak{n}}^{\mathfrak{b}}(\lambda)}_{\xi}$
is locally finite. Motivated by the results from Subsection~\ref{s6.1} 
it is natural to expect that the algebra describing  $\overline{\mathfrak{W}_{\mathfrak{n}}^{\mathfrak{b}}(\lambda)}_{\xi}$
is Koszul. From the proof of Corollary~\ref{cor7015} it is easy to see 
this algebra is more complicated than the tensor products of algebras from
Subsection~\ref{s6.1} (such tensor products are obviously Koszul).
If one makes a parallel with simple finite dimensional Lie algebras, 
then our expectation of Koszulity for this algebra is similar 
to Alexandru conjecture (for thick category $\mathcal{O}$), see
\cite{Gai}.

\subsection{Some solvable subquotients of the Virasoro algebra}\label{s6.4}

Another way to generalize the results of Subsection~\ref{s6.1} is to consider
certain subquotients of the Virasoro algebra. For $n=0,1,2,\dots$ let
$\mathfrak{v}_n$ denote the Lie algebra with the basis 
$\{e_i:i=n,n+1,\dots\}$,
and the Lie bracket given by $[e_i,e_j]=(j-i)e_{i+j}$. For $n>0$ the
algebra $\mathfrak{v}_n$ is quasi-nilpotent. For $k\geq n$ the algebra
$\mathfrak{v}_k$ is an ideal of $\mathfrak{v}_n$. The quotient
$\mathfrak{v}_n/\mathfrak{v}_k$ is always solvable and, moreover, nilpotent
if $n>0$. In particular, it is easy to see that 
$(\mathfrak{v}_n,\mathfrak{v}_k)$ is a Whittaker pair for all $k>n$
and $(\mathfrak{v}_n/\mathfrak{v}_k,\mathfrak{v}_m/\mathfrak{v}_k)$
is a Whittaker pair for all $k>m>n$. For $n=0$, $m=1$ and $k=2$ one
obtains the algebra considered in Subsection~\ref{s6.1}.
For $n=1$, $m=3$ and $k=4$ one obtains the algebra considered 
in Subsection~\ref{s6.2}.

For all Whittaker pairs $(\mathfrak{g},\mathfrak{n})$ of the 
form $(\mathfrak{v}_0,\mathfrak{v}_1)$
and $(\mathfrak{v}_0/\mathfrak{v}_k,\mathfrak{v}_1/\mathfrak{v}_k)$
the module $M_{\lambda}$ is isomorphic to $\mathbb{C}[e_0]$ as
a $\mathbb{C}[e_0]$-module. It is simple if and only if $\lambda\neq 0$
(the ``if'' part follows from Proposition~\ref{prop6001} and the 
``only if'' part is obvious).
If $\lambda=0$, then $M_{\lambda}$ is free of rank one over its endomorphism
algebra $\mathrm{End}_{U(\mathfrak{g})}(M_{\lambda})\cong
\mathbb{C}[e_0]$ and simple quotients of $M_{\lambda}$ are all
one-dimensional and have the form $L_{\mu}$, where $\mu\in\mathbb{C}$,
$\mathfrak{n}L_{\mu}=0$ and $e_0v=\mu v$ for all $v\in L_{\mu}$.
Similarly to Corollary~\ref{cor7015} one can show that 
indecomposable blocks of the category
$\overline{\mathfrak{W}_{\mathfrak{n}}^{\mathfrak{b}}}$
can be described as module categories over 
(completions of) some locally finite quiver algebras with relations.
We believe that these algebras are Koszul.

\vspace{2cm}

\noindent
Punita Batra, Harish-Chandra Research Institute,
Chhatnag Road,  Jhusi, Allahabad, 211 019, INDIA,
{\tt batra\symbol{64}mri.ernet.in}
\vspace{0.8cm}

\noindent
Volodymyr Mazorchuk, Department of Mathematics, Uppsala University,
Box 480, 751 06, Uppsala, SWEDEN, {\tt mazor\symbol{64}math.uu.se}\\
http://www.math.uu.se/$\tilde{\hspace{1mm}}$mazor/.
\vspace{0.8cm}


\begin{thebibliography}{9999}
\bibitem[AP]{AP} D.~Arnal, G.~Pinczon; On algebraically irreducible
representations of the Lie algebra ${\mathfrak sl}(2)$.  
J. Mathematical Phys.  {\bf 15}  (1974), 350--359.
\bibitem[Ba]{Ba} E.~Backelin; Representation of the category 
$\mathcal O$ in Whittaker categories.  Internat. Math. Res. Notices  
1997,  no. {\bf 4}, 153--172.
\bibitem[Bav]{Bav} V.~Bavula; The extension group of the simple modules 
over the first Weyl algebra.  Bull. London Math. Soc.  {\bf 32}  
(2000),  no. 2, 182--190.
\bibitem[BO]{BO} G.~Benkart, M.~Ondrus; Whittaker modules for 
generalized Weyl algebras. Represent. Theory {\bf 13} (2009), 141--164.
\bibitem[Bl]{Bl} R.~Block; The irreducible representations of the Lie 
algebra $\mathfrak{sl}(2)$ and of the Weyl algebra.  Adv. in Math.  
{\bf 39}  (1981), no. 1, 69--110.
\bibitem[CE]{CE} H.~Cartan, S.~Eilenberg; Homological algebra. With an 
appendix by David A. Buchsbaum. Reprint of the 1956 original.
Princeton Landmarks in Mathematics. Princeton University Press, 
Princeton, NJ, 1999.
\bibitem[Ch]{Ch} K.~Christodoulopoulou; Whittaker modules for 
Heisenberg algebras and imaginary Whittaker modules for affine 
Lie algebras.  J. Algebra  {\bf 320}  (2008),  no. 7, 2871--2890.
\bibitem[CM]{CM} C.~Conley, C.~Martin, Annihilators of tensor density 
modules. J. Algebra {\bf 312} (2007), no. 1, 495--526. 
\bibitem[Di]{Di} J.~Dixmier, Enveloping algebras. Revised reprint of 
the 1977 translation. Graduate Studies in Mathematics, {\bf 11}. 
American Mathematical Society, Providence, RI, 1996.
\bibitem[DFO]{DOF} Yu.~Drozd, V.~Futorny, S.~Ovsienko; 
Harish-Chandra subalgebras and Gelfand-Zetlin modules.  
Finite-dimensional algebras and related topics (Ottawa, ON, 1992),  
79--93, NATO Adv. Sci. Inst. Ser. C Math. Phys. Sci., {\bf 424}, 
Kluwer Acad. Publ., Dordrecht, 1994.
\bibitem[Ga]{Ga} P.~Gabriel; Indecomposable representations. II.  
Symposia Mathematica, Vol. XI (Convegno di Algebra Commutativa, INDAM, 
Rome, 1971),  pp. 81--104. Academic Press, London, 1973.
\bibitem[Gai]{Gai} P.-Y.~Gaillard; Statement of the Alexandru 
Conjecture, Preprint arXiv:math/0003070.
\bibitem[KM]{KM} O.~Khomenko, V.~Mazorchuk; Structure of modules 
induced from simple modules with minimal annihilator.  Canad. J. Math.  
{\bf 56}  (2004),  no. 2, 293--309. 
\bibitem[Ko]{Ko} B.~Kostant; On Whittaker vectors and representation 
theory. Invent. Math. {\bf 48} (1978), no. 2, 101--184. 
\bibitem[LW]{LW} J.~Li, B.~Wang; Whittaker Modules For The $W$-algebra 
$W(2,2)$. Preprint arXiv:0902.1592.
\bibitem[LWZ]{LWZ} D.~Liu, Y.~Wu, L.~Zhu; Whittaker Modules for the 
twisted Heisenberg-Virasoro Algebra. Preprint arXiv:0902.4074.
\bibitem[Ma]{Ma} O.~Mathieu; Classification of irreducible weight 
modules.  Ann. Inst. Fourier (Grenoble)  {\bf 50}  (2000),  no. 2, 537--592.
\bibitem[MOS]{MOS} V.~Mazorchuk, S.~Ovsienko, C.~Stroppel;
Quadratic duals, Koszul dual functors, and applications.  
Trans. Amer. Math. Soc.  {\bf 361}  (2009),  no. 3, 1129--1172. 
\bibitem[MD1]{MD1} E.~McDowell; On modules induced from Whittaker 
modules.  J. Algebra  {\bf 96}  (1985),  no. 1, 161--177. 
\bibitem[MD2]{MD2} E.~McDowell; A module induced from a Whittaker 
module.  Proc. Amer. Math. Soc.  {\bf 118}  (1993),  no. 2, 349--354.
\bibitem[MS1]{MS1} D.~Mili{\v c}i{\'c}, W.~Soergel; The composition 
series of modules induced from Whittaker modules.  Comment. Math. Helv.  
{\bf 72} (1997),  no. 4, 503--520.
\bibitem[MS2]{MS2} D.~Mili{\v c}i{\'c}, W.~Soergel; Twisted 
Harish-Chandra sheaves and Whittaker modules: The non-degenerate case.
Preprint, available from 
http://home.mathematik.uni-freiburg.de/soergel/
\bibitem[MP]{MP} R.~Moody, A.~Pianzola; Lie algebras with triangular
decompositions. Canadian Mathematical Society Series of Monographs 
and Advanced Texts. A Wiley-Interscience Publication. John Wiley \& 
Sons, Inc., New York, 1995. 
\bibitem[On1]{On1} M.~Ondrus; Whittaker modules for 
$U\sb q({\mathfrak{sl}}\sb 2)$.  J. Algebra  {\bf 289}  
(2005),  no. 1, 192--213.
\bibitem[On2]{On2} M.~Ondrus; Tensor products and Whittaker vectors 
for quantum groups.  Comm. Algebra  {\bf 35}  (2007),  no. 8, 2506--2523.
\bibitem[OW]{OW} M.~Ondrus, E.~Wiesner; Whittaker Modules for the 
Virasoro Algebra. J. Algebra Appl. {\bf 8} (2009), no. 3, 363-377.
\bibitem[PZ]{PZ} I.~Penkov, G.~Zuckerman; Generalized Harish-Chandra 
modules: a new direction in the structure theory of representations.  
Acta Appl. Math.  {\bf 81}  (2004),  no. 1-3, 311--326.
\bibitem[Ru]{Ru} A.~Rudakov; Irreducible representations of 
infinite-dimensional Lie algebras of Cartan type. Izv. Akad. Nauk 
SSSR Ser. Mat.  {\bf 38}  (1974), 835--866. 
\bibitem[TZ]{TZ} S.~Tan, X.~Zhang; Whittaker modules for 
the Schr{\"o}dinger-Virasoro algebra. Preprint arXiv:0812.3245.
\bibitem[Ta]{Ta} X.~Tang; On Whittaker modules over a class of algebras 
similar to $U({\mathfrak sl}\sb 2)$. Front. Math. China {\bf 2} (2007), 
no. 1, 127--142. 
\bibitem[Vo]{Vo} D.~Vogan, Jr.; Representations of real reductive Lie 
groups. Progress in Mathematics, {\bf 15}. Birkh{\"a}user, Boston, Mass., 
1981. 
\bibitem[Wa]{Wa} B.~Wang; Whittaker Modules for Graded Lie Algebras.
Preprint arXiv:0902.3801.
\bibitem[WZ]{WZ} B.~Wang, X.~Zhu; Whittaker modules for a Lie algebra 
of Block type. Preprint arXiv:0907.0773.
\end{thebibliography}
\end{document}